\documentclass[11pt, reqno]{amsart}
\usepackage{amsmath, amsfonts, amsthm, amssymb, multicol, mathtools, dsfont,verbatim}
\usepackage{graphicx}
\usepackage{float}
\usepackage{scalerel,stackengine, subcaption}
\usepackage[usenames,dvipsnames,x11names]{xcolor}
\usepackage{enumitem}
\usepackage{pgfplots}
\usepgfplotslibrary{fillbetween}
\pgfplotsset{width=10cm,compat=1.9}
\usepackage[square,sort,comma,numbers]{natbib}
\usepackage[pdfencoding=auto, psdextra]{hyperref}
\makeatletter
\def\@setauthors{%
\begingroup
\def\thanks{\protect\thanks@warning}%
\trivlist
\centering\footnotesize \@topsep30\p@\relax
\advance\@topsep by -\baselineskip
  \item\relax
\author@andify\authors
\def\\{\protect\linebreak}

\normalsize\lowercase{\authors}%

\ifx\@empty\contribs
\else
,\penalty-3 \space \@setcontribs
\@closetoccontribs
\fi
\endtrivlist
\endgroup
}
\def\@settitle{\begin{center}
\LARGE\lowercase{\@title}
  \end{center}%
}
\makeatother
\newcommand{\authoremail}[1]{\email{\href{mailto:#1}{\color{lightblue}{#1}}}}
\newcommand{\authoraddress}[1]{\address{\normalfont{#1}}}

\usepackage{fancyhdr}
\numberwithin{equation}{section}

\newtheorem{thm}{Theorem}[section]

\newtheorem{cor}[thm]{Corollary}

\newtheorem{ques}[thm]{Question}

\DeclareMathOperator{\dist}{dist}

\newtheorem{theorem}{Theorem}[section]

\theoremstyle{definition}

\newtheorem{example}[theorem]{Example}

\theoremstyle{remark}

\renewcommand{\epsilon}{\varepsilon}

\newcommand{\R}{\mathbb{R}}

\renewcommand{\ge}{\csname geqslant\endcsname}
\renewcommand{\le}{\csname leqslant\endcsname}
\renewcommand{\geq}{\csname geqslant\endcsname}
\renewcommand{\leq}{\csname leqslant\endcsname}

\newcommand{\ubd}{\overline{\dim}_{\textup{B}}\,}
\newcommand{\ubdzero}{\overline{\dim}_{\textup{B}}^0\,}
\newcommand{\lbd}{\underline{\dim}_{\textup{B}}\,}

\newcommand{\ad}{\dim_{\mathrm{A}}}
\newcommand{\qad}{\dim_{\mathrm{qA}}}
\newcommand{\as}{\dim^\theta_{\mathrm{A}}}

\newcommand{\hd}{\dim_{\mathrm{H}}}

\newcommand{\rn}{\mathbb{R}^d}

\newcommand{\diam}{\textnormal{diam}}

\renewcommand{\epsilon}{\varepsilon}

\newcommand{\rd}{\mathbb{R}^d}

\newcommand{\gdk}{G(d,k)}

\newcommand{\gamdk}{\gamma_{d,k}}

\newcommand{\rk}{\mathbb{R}^k}

\newcommand{\Q}{\mathbb{Q}}
\newcommand{\N}{\mathbb{N}}
\renewcommand{\L}{\mathcal{L}}

\newcommand{\spt}{\text{spt}\,}

\DeclareMathOperator*{\esssup}{ess\,sup}

\makeatletter
\DeclareRobustCommand\widecheck[1]{{\mathpalette\@widecheck{#1}}}
\def\@widecheck#1#2{%
\setbox\z@\hbox{\m@th$#1#2$}%
\setbox\tw@\hbox{\m@th$#1%
\widehat{%
\vrule\@width\z@\@height\ht\z@
\vrule\@height\z@\@width\wd\z@}$}%
\dp\tw@-\ht\z@
\@tempdima\ht\z@ \advance\@tempdima2\ht\tw@ \divide\@tempdima\thr@@
\setbox\tw@\hbox{%
\raise\@tempdima\hbox{\scalebox{1}[-1]{\lower\@tempdima\box
\tw@}}}%
{\ooalign{\box\tw@ \cr \box\z@}}}
\makeatother

\stackMath
\newcommand\reallywidehat[1]{%
\savestack{\tmpbox}{\stretchto{%
\scaleto{%
\scalerel*[\widthof{\ensuremath{#1}}]{\kern.1pt\mathchar"0362\kern.1pt}%
{\rule{0ex}{\textheight}}
}{\textheight}%
}{2.4ex}}%
\stackon[-6.9pt]{#1}{\tmpbox}%
}
\definecolor{lightblue}{HTML}{2B77A4}
\colorlet{plotblue}{LightSkyBlue3!80}
\definecolor{darkred}{HTML}{9E0D0D}
\definecolor{darkyellow}{HTML}{b3b300}
\definecolor{darkorange}{HTML}{D86129}
\hypersetup{
colorlinks=true,
linkcolor=darkred,
urlcolor=darkred,
citecolor=lightblue
}
\title{On the Marstrand projection theorem for the Assouad spectrum}

\author{Kenneth J. Falconer}
\authoraddress{K. J. Falconer, University of St Andrews, Scotland}
\authoremail{kjf@st-andrews.ac.uk}

\author{Jonathan M. Fraser}
\authoraddress{J. M. Fraser, University of St Andrews, Scotland}
\authoremail{jmf32@st-andrews.ac.uk}
\thanks{JMF was financially supported by  a  \emph{Leverhulme Trust Research Project Grant} (RPG-2023-281) and an  \emph{EPSRC Open Fellowship} (EP/Z533440/1).}

\author{Antti K\"aenm\"aki}
\authoraddress{A. K\"aenm\"aki, University of Eastern Finland, Finland} 
\authoremail{antti@kaenmaki.net}
\thanks{ }

\begin{document}
\thispagestyle{empty}

\begin{abstract}
  Marstrand's projection theorem states that the Hausdorff dimension of the orthogonal projection of a Borel set in the plane onto lines is constant almost surely. This property extends to other notions of dimension, such as box  and packing dimensions, but does not hold for the Assouad dimension. In this paper, we show that Marstrand's projection theorem also fails for the quasi-Assouad dimension and the Assouad spectrum, which interpolates between the upper box  and quasi-Assouad dimensions. Additionally, we establish an almost sure lower bound for the Assouad spectrum of the projections using capacity-theoretic dimension profiles, and an almost sure upper bound for projections of bounded planar sets via an incidence geometry inspired  tube-counting argument. As an application, for a parametrised family of homogeneous self-similar sets, we obtain an almost sure upper bound for the Assouad spectrum which beats the trivial upper bound coming from the upper box dimension. 
\end{abstract}

\maketitle

\tableofcontents

\section{Introduction} \label{intro}

\subsection{Dimension theory}

Both the Hausdorff and box dimensions measure the size of a fractal set via covers of the whole set. In many settings---such as embedding theory \cite{robinson}---finer local information is needed, and the \emph{Assouad dimension} captures it. For $X \subseteq \rd$, it is defined by
\begin{align*}
  \dim_\text{A} X  \ = \ \inf \bigg\{ \ \alpha &: \text{   there exists  a constant  $C >0$ such that} \\
  &\text{$ N_{r} \big( X \cap B(x,R) \big) \ \leq \ C \bigg(\frac{R}{r}\bigg)^\alpha$ }
  \text{for all $0< r<R$ and $x \in X$} \bigg\},
\end{align*}
where $N_r(A)$ is the smallest number of closed balls of radius $r>0$ needed to cover the bounded set $A$. Fraser and Yu \cite{assouadspectrum} introduced the Assouad spectrum to interpolate between the upper box dimension and the Assouad dimension by using a parameter $\theta \in (0,1)$ to fix the relationship between the two scales $r<R$ in the definition of the Assouad dimension. More precisely, the \emph{Assouad spectrum} of $X \subseteq \rd$ at $\theta \in (0,1)$ is defined by
\begin{align*}
  \as X  \ = \ \inf \bigg\{ \ \alpha &: \text{   there exists  a constant  $C >0$ such that} \\
  &\text{$ N_{R^{1/\theta}} \big( X \cap B(x,R) \big) \ \leq \ C \bigg(\frac{R}{R^{1/\theta}}\bigg)^\alpha$ }
  \text{for all $0<  R<1$ and $x \in X$} \bigg\}.
\end{align*}
The Assouad spectrum, which is bi-Lipschitz invariant, may be viewed as the map $\theta \mapsto \as X$. The \emph{quasi-Assouad dimension} is defined by
\begin{equation*}
  \qad X \ = \ \lim_{\theta \uparrow 1} \as X.
\end{equation*}
The quasi-Assouad dimension was introduced by L\"u and Xi in \cite{quasiassouad}; the characterisation above in terms of the Assouad spectrum was proved in \cite{canadian}. Finally, the \emph{upper box  dimension} of a bounded set $X$ is
\begin{equation*}
  \ubd X \ = \ \limsup_{r \downarrow 0} \frac{\log N_r(X)}{-\log r}.
\end{equation*}
The lower box  dimension $\lbd X$ is defined analogously with $\liminf_{r \downarrow 0}$. The key distinction between the Assouad dimension and the box dimension is that the Assouad dimension counts covers within a single ball rather than of the whole set, with the covering number normalised accordingly. Generally, for $\theta \in (0,1)$ and sets $X$,
\begin{equation*}
  \ubd X \ \leq \ \as X \ \leq \ \qad X \ \leq \ \ad X
\end{equation*}
and
\begin{equation*}
  \as X \ \leq \ \min\bigg\{\frac{\ubd X}{1-\theta} , \ \qad X\bigg\},
\end{equation*}
where $X$ is assumed to be bounded for the inequalities involving the upper box dimension.
In particular,
\begin{equation*}
  \ubd X \ = \ \lim_{\theta \downarrow 0} \as X,
\end{equation*}
which may be viewed as continuity at the endpoint $\theta=0$. The possible functions realised as Assouad spectra were completely classified in \cite{specclass}.

The definition of the upper box dimension can be extended to an unbounded set $X$ simply by taking $\sup_{n \in \N} \ubd(X \cap B(0,n))$. Wang and Li \cite{wangli} recently proposed another such extension, namely
\begin{equation*}
  \ubdzero X \ = \ \lim_{\theta \downarrow 0} \as X.
\end{equation*}
They showed that $\ubd X \leq \ubdzero X$, where the inequality can be strict for unbounded sets. Since the modified dimension
\begin{equation*}
  \inf\biggl\{ \sup_{n \in \N} \ubdzero X_n : X \subseteq \bigcup_n X_n \biggr\}
\end{equation*}
equals the packing dimension of $X$, they suggested that $\ubdzero$ would be a more natural definition for the upper box dimension of an unbounded set.

The Assouad spectrum has diverse applications, including $L^p \to L^q$ bounds for spherical maximal functions \cite{anderson, roos, beltran}, weak embeddability problems \cite{stathis}, certain H\"older regularity questions \cite{holder}, quasiconformal mapping problems \cite{stathistyson,stathisquasi}, and connections to the Sullivan dictionary in conformal dynamics \cite{bullams}.

\subsection{Dimensions of orthogonal projections}

For integers $1\leq k<d$, we write $G(d,k)$ for the Grassmannian manifold of $k$-dimensional subspaces $V\subseteq\rd$.  This space may be equipped with the invariant Borel probability measure $\gamdk$ obtained from the Haar measure on the topological group of rotations around the origin. We write  $P_{V}\colon \rd \to V$ for the orthogonal projection onto the subspace $V\in\gdk$, which we identify with $\rk$.

A central question in fractal geometry and geometric measure theory is how the dimensions of a set $X \subseteq \rd$ compare to those of its orthogonal projections $P_V X$. Since the Hausdorff dimension cannot increase under Lipschitz maps, we have
\begin{equation*}
  \hd P_V X \ \leq \ \min\{k, \hd X\}
\end{equation*}
for \emph{all} $V \in G(d,k)$. The Assouad and quasi-Assouad dimensions and the Assouad spectrum can increase under Lipschitz maps, so the simple bound above does not hold for these notions. Marstrand's projection theorem is one of the most well-known and influential results in fractal geometry and geometric measure theory. Marstrand \cite{Mar54} proved it in the plane, Kaufman \cite{Kau68} gave a potential-theoretic proof, and Mattila \cite{Mat75} extended it to higher dimensions, showing that for Borel sets $X\subseteq\rd$ and $\gamdk$-almost all $V\in\gdk$,
\begin{equation*}
  \hd P_V X \ = \ \min\{k, \hd X \}.
\end{equation*}
This theorem stimulated substantial activity in fractal geometry, geometric measure theory, and related areas; see \cite{FFJ15, Mat14,Shm15, falconer25}.

Marstrand's theorem invites the question of analogous almost-sure results for other notions of dimension. For box dimensions, the answer is affirmative: orthogonal projections have an almost-surely constant box dimension, but the value is given by a more intricate ``dimension profile'' rather than simply $\min\{k,\hd X\}$. These results were first established by Falconer and Howroyd \cite{falconerhowroyd,falconerhowroyd2,howroyd}; see also \cite{falconerprofile, falconerprofile2}.

Dimension profiles can be defined in terms of capacities with respect to certain kernels \cite{falconerprofile}. For $0 \leq s \leq d$ and $r>0$ we define the kernel
\begin{equation} \label{ker}
  \phi_r^s(x) \ = \ \min\bigg\{ 1, \bigg(\frac{r}{|x|}\bigg)^s\bigg\}, \qquad x\in \rn.
\end{equation}
For a non-empty compact $X \subseteq \rn$, the \emph{capacity} $C_r^s(X)$ with respect to this kernel is given by
\begin{equation} \label{eq:capacity-def}
  C_r^s(X)^{-1}\ = \ \inf_{\mu \in {\mathcal M}(X)}\iint \phi_r^s(x-y)\,\mathrm{d}\mu(x)\,\mathrm{d}\mu(y)
\end{equation}
where ${\mathcal M}(X)$ denotes the set of Borel probability measures supported by $X$. The double integral is the {\it energy} of $\mu$ relative to the kernel, analogous to the $s$-energy $I_s(\mu)$ in potential-theoretic approaches to Hausdorff dimension; see e.g.~\cite[Theorem 4.13]{Fal03}. The capacity of a bounded set is taken to be that of its closure.

For bounded $X \subseteq \rn$ and $s>0$, we define the {\it lower} and {\it upper box dimension profiles} of $X$ by
\begin{equation*}
  \underline{\mbox{\rm dim}}_{\rm B}^s X \ = \ \varliminf_{r\downarrow 0} \frac{\log C_r^s(X)}{-\log r} \quad\text{and}\quad \overline{\mbox{\rm dim}}_{\rm B}^s X \ = \ \varlimsup_{r\downarrow 0} \frac{\log C_r^s(X)}{-\log r}.
\end{equation*}
In particular, by \cite[Corollary 2.5]{falconerprofile} if $s\geq d$ then
\begin{equation*}
  \underline{\mbox{\rm dim}}_{\rm B}^s X \ = \ \lbd X \quad\text{and}\quad \overline{\mbox{\rm dim}}_{\rm B}^s X \ = \ \ubd X,
\end{equation*}
but for $s < d$ the dimension profiles give the almost-sure dimensions of projections of sets as well as the size of the exceptional set.

\begin{thm}{\rm \cite[Theorems 1.1 and 1.2]{falconerprofile}}\label{mainA}

  \noindent $(i)$ Let $1\leq k<d$ be an integer.  For almost all $V \in G(d,k)$, if $X\subseteq \rn$ is bounded
  \begin{equation*}
    \lbd P_V X \ = \ \underline{\mbox{\rm dim}}_{\rm B}^k X \quad\text{and}\quad \ubd P_V X \ = \ \overline{\mbox{\rm dim}}_{\rm B}^k X.
  \end{equation*}
  \noindent $(ii)$ For $0<s<k$, if $X\subseteq \rn$ is bounded
  \begin{align*}
    \hd \{ V \in G(d,k) : \lbd P_V X < \underline{\mbox{\rm dim}}_{\rm B}^s X\} \ &\leq \ k(d-k)-(k-s), \\
    \hd \{ V \in G(d,k) : \ubd P_V X < \overline{\mbox{\rm dim}}_{\rm B}^s X\} \ &\leq \ k(d-k)-(k-s).
  \end{align*}
\end{thm}

Although the values of $\lbd P_V X$ and $ \ubd P_V X$ are constant for almost all $V\in G(d,k)$, this constant can take any value in the range
\begin{equation} \label{boxdimsbounds}
  \frac{\ubd X}{1+(\frac{1}{k}-\frac{1}{d})\ubd X} \ \leq \ \overline{\mbox{\rm dim}}_{\rm B}^k X \ \leq \ \min \{k, \ubd X\},
\end{equation}
with analogous inequalities for lower box dimension. These inequalities were established in \cite{falconerhowroyd, falconerhowroyd2, howroyd} using dimension profiles directly, with examples showing them to be best possible in \cite{falconerhowroyd}. A simpler approach using capacities was recently given in \cite{falconerprofile}.

There is no Marstrand theorem for the Assouad dimension: the Assouad dimension of the orthogonal projection of a compact set in $\rd$ need not take a constant value almost surely. This was established by Fraser and Orponen \cite{FO17}, and further extreme behaviour was exhibited in \cite{antti}, including examples where $\ad P_V X$ takes a different value for each distinct $V$. Despite the absence of a Marstrand theorem for the Assouad dimension, for $\gamdk$-almost all $V \in G(d,k)$,
\begin{equation} \label{assthm}
  \ad P_V X \ \geq \ \min\{k, \ad X\},
\end{equation}
providing the expected lower bound almost surely. This was proved in the plane in \cite{FO17} and extended to arbitrary dimensions in \cite{fraserisrael}. Orponen \cite{orponenassouad} refined the planar case by showing that the set of directions $V$ violating \eqref{assthm} has not only measure zero but also Hausdorff dimension zero. Recently, Wu \cite{wu} further strengthened this, demonstrating that the set of exceptional directions is countable. In this paper, we establish that the Marstrand projection theorem also fails for the Assouad spectrum and quasi-Assouad dimension, and we provide almost-sure lower and upper bounds for the spectrum of projections.

\section{Main results} \label{sec:main-results}

\subsection{Almost sure lower bounds for the Assouad spectrum of projections}

We first establish almost-sure lower bounds for the Assouad spectrum of projections in terms of dimension profiles. For $s\leq d$ and $ r > 0$, let the kernel $\phi_r^s : \R^d \to \R$ be as in \eqref{ker} and the capacity $C_r^s$ as in \eqref{eq:capacity-def}. Standard potential theory ensures the existence of a measure $\mu \in \mathcal{M}(X)$, called an \emph{equilibrium measure} associated with the capacity $C_r^s(X)$, such that
\begin{equation*}
  C_r^s(X)^{-1}\ = \ \iint \phi_r^s(x-y)\,\mathrm{d}\mu(x)\,\mathrm{d}\mu(y).
\end{equation*}
Combining the covering definition of the Assouad spectrum with \cite[Corollary 2.4]{falconerprofile}, and absorbing the logarithmic factor in the case $s = d$ into the exponent, the Assouad spectrum can equivalently be defined as
\begin{equation} \label{eq:assouad-spectrum-characterization}
\begin{split}
  \as X  \ = \ \inf \bigg\{ \ \alpha &: \text{   there exists  a constant  $C >0$ such that} \\
  &C_{R^{1/\theta}}^d(X \cap B(x, R)) \leq C \biggl(\frac{R}{R^{1/\theta}}\biggr)^\alpha
  \text{for all $0<  R<1$ and $x \in X$} \bigg\}.
\end{split}
\end{equation}
Define the \emph{$k$-dimensional Assouad spectrum profile of $X$} by
\begin{align*}
  \dim_{\textup{A}}^{\theta,k} X \ = \ \inf \bigg\{ \ \alpha &: \text{   there exists  a constant  $C >0$ such that} \\
  &C_{R^{1/\theta}}^k(X \cap B(x, R)) \leq C \biggl(\frac{R}{R^{1/\theta}}\biggr)^\alpha
  \text{for all $0<  R<1$ and $x \in X$} \bigg\}.
\end{align*}
Analogous to the box dimension case, these profiles can be used to study projections onto $k$-dimensional subspaces. In particular, we obtain the following almost-sure lower bound.

\begin{thm}\label{main2}
  Let $X \subset \R^d$, $\theta \in (0,1)$, and $k < d$ be an integer. Then
  \begin{equation*}
    \as P_V X \ \geq \ \dim_{\textup{A}}^{\theta,k} X
  \end{equation*}
  for $\gamma_{d,k}$-almost all $V \in G(d,k)$.
\end{thm}

Theorem \ref{main2} is most useful when the profile $\dim_{\textup{A}}^{\theta,k} X$ can be bounded from below by a more familiar invariant. Corollary \ref{cor:assouad-spectrum-profile-bound} gives such an explicit bound in terms of $\as X$, at the cost of changing the spectrum parameter from $\theta$ to $\theta'$.

\begin{cor} \label{cor:assouad-spectrum-profile-bound}
  Let $X \subset \R^d$, $\theta \in (0,1)$, and $1 \leq k < d$ be an integer. Let $\theta' \in (0,1)$ be defined by
  \begin{equation*}
    \frac{1}{\theta'} -1 \ = \ \biggl(\frac{1}{\theta}-1\biggr)\biggl(1+\biggl(\frac{1}{k}-\frac{1}{d}\biggr)\as X\biggr).
  \end{equation*}
  Then
  \begin{equation*}
    \dim_{\textup{A}}^{\theta'} P_V X \ \geq \ \frac{\as X}{1+(\frac{1}{k}-\frac{1}{d})\as X}
  \end{equation*}
  for $\gamma_{d,k}$-almost all $V \in G(d,k)$.
\end{cor}

The proof of Theorem \ref{main2}, and the derivation of Corollary \ref{cor:assouad-spectrum-profile-bound} from it, are given in Section \ref{sec:almost-sure-lower-bound}. The endpoints $\theta = 0$ and $\theta = 1$ are not covered by Corollary \ref{cor:assouad-spectrum-profile-bound}, but both are recovered as limiting cases. If $X$ is bounded and $\theta \downarrow 0$, then $\as X \to \ubd X$ and the defining relation for $\theta'$ gives $\theta' \downarrow 0$. Taking limits in Corollary \ref{cor:assouad-spectrum-profile-bound} along a countable sequence of full-measure sets gives
\begin{equation*}
  \ubd P_V X \ \geq \ \frac{\ubd X}{1+(\frac{1}{k}-\frac{1}{d})\ubd X}
\end{equation*}
for $\gamma_{d,k}$-almost all $V \in G(d,k)$. Together with Theorem \ref{mainA}(i), this is the lower estimate in \eqref{boxdimsbounds}. The H\"older estimate used in the proof is the capacity interpolation that makes this limiting agreement with the box dimension profile possible. At the other endpoint, $\theta \uparrow 1$ implies $\theta' \uparrow 1$, although $\theta' \neq \theta$ in general for fixed $\theta < 1$. Taking this limit in Theorem \ref{main2} and Corollary \ref{cor:assouad-spectrum-profile-bound} yields corresponding lower bounds for the quasi-Assouad dimension. We define the \emph{$k$-dimensional quasi-Assouad dimension profile of $X$} by
\begin{equation*}
  \dim_{\textup{qA}}^{k} X \ = \ \limsup_{\theta \uparrow 1} \dim_{\textup{A}}^{\theta,k} X.
\end{equation*}
This immediately yields the following corollary.

\begin{cor}\label{main2cor}
  Let $X \subset \R^d$ and $1 \leq k < d$ be an integer. Then
  \begin{equation*}
    \qad P_V X \ \geq \ \dim_{\textup{qA}}^{k} X \ \geq \ \frac{\qad X}{1+(\frac{1}{k}-\frac{1}{d})\qad X}
  \end{equation*}
  for $\gamma_{d,k}$-almost all $V \in G(d,k)$.
\end{cor}

\subsection{Tube-counting upper bounds in the plane}

We now turn from general lower bounds to upper bounds for projections of bounded planar sets. The next result gives an almost-sure upper bound for $\as P_V X$, improving on the trivial estimate
\begin{equation*}
  \as P_V X \ \leq \ \frac{ \ubd P_V X }{1-\theta} \ \leq \ \frac{ \ubd X }{1-\theta}
\end{equation*}
whenever
\begin{equation*}
  \as X \ < \ \frac{ \ubd X }{1-\theta}.
\end{equation*}
This result is proved by a tube-counting approach inspired by incidence geometry.  Consider the projections of a bounded set in the plane  and suppose many of the projections have large Assouad spectrum.  Then, for each such projection, one can find small balls which are highly concentrated.  By considering incidences between the  tubes given by the preimages of these highly concentrated balls under the projection map,   one can select well-separated points from   the  tubes   and combine them into a large well-separated subset of $X$, forcing $\ubd X$ to be large. Thus the set $X$ itself must already have sufficiently large dimension,  ruling out many potential bounded counterexamples. In particular, this shows that the behaviour of the Assouad spectrum under projections is not `maximally wild' as in the case of the Assouad dimension. 

\begin{thm} \label{thm:tube-lower-spectrum}
  Let $X \subset \R^2$ be bounded and $\theta \in (0,1)$. Then
  \begin{equation*}
    \as P_V X \ \leq \ \frac{1}{2} \biggl( \frac{\ubd X}{1-\theta} + \as X \biggr)
  \end{equation*}
  for $\gamma_{2,1}$-almost all $V \in G(2,1)$.
\end{thm}

See Section \ref{sec:tube-counting-obstructions} for the proof. For example, this upper bound shows that if $X$ is bounded, $\theta \in (0,\frac23)$, and $\as X < \varepsilon$, then $\as P_V X \leq 2\varepsilon$ for $\gamma_{2,1}$-almost all $V \in G(2,1)$. Thus, for bounded sets, the Assouad spectrum of projections cannot be large almost everywhere if $X$ itself has small dimension. The boundedness hypothesis enters the proof in two ways. The easy case $\esssup_{V \in G(2,1)} \as P_V X \leq \as X$ is closed by the standard inequality $\as X \leq \ubd X/(1-\theta)$, which we use only for bounded sets. More importantly, after extracting many $2r$-separated points from the exceptional tubes, boundedness allows one to convert these points into a lower bound for the global covering number $N_r(X)$ and hence for $\ubd X$; for unbounded sets, points selected from different tubes may lie arbitrarily far apart in the ambient plane, so this final step need not yield a contradiction. In Theorems \ref{main1}, \ref{main1+}, and \ref{maintheta} we construct unbounded examples where this phenomenon occurs.

The following result provides a trade-off between the upper bound and the Hausdorff dimension of the exceptional set.

\begin{thm} \label{thm:tube-lower-spectrum-alpha}
  Let $X \subset \R^2$ be bounded, $\theta \in (0,1)$, and $\alpha \in [0,1)$. Then
  \begin{equation*}
    \hd \biggl\{ \ V \in G(2,1) : \ \as P_V X > \frac{1}{1+\alpha} \biggl( \frac{\ubd X}{1-\theta} + \alpha \as X \biggr)\biggr\} \ \leq \ \alpha.
  \end{equation*}
\end{thm}

See Section \ref{sec:tube-counting-obstructions-alpha} for the proof. Inserting $\alpha = 1$ formally into the threshold above yields the upper bound from Theorem \ref{thm:tube-lower-spectrum}, but the conclusion guarantees only that the exceptional set has Hausdorff dimension at most $1$. Thus, Theorem \ref{thm:tube-lower-spectrum-alpha} does not recover Theorem \ref{thm:tube-lower-spectrum} in the limit as $\alpha \uparrow 1$. Nevertheless, the limit $\alpha \downarrow 0$ yields an endpoint statement. If $X \subset \R^2$ is bounded with $\as X < \ubd X/(1-\theta)$, then
\begin{equation} \label{eq:zerodim}
  \hd \biggl\{ \ V \in G(2,1) : \ \as P_V X = \frac{\ubd X}{1-\theta} \biggr\} \ = \ 0.
\end{equation}
Indeed, for each $\alpha \in (0,1)$, the number
\begin{equation*}
  u_\alpha \ = \ \frac{1}{1+\alpha} \biggl( \frac{\ubd X}{1-\theta} + \alpha \as X \biggr)
\end{equation*}
satisfies $u_\alpha < \ubd X/(1-\theta)$, and therefore $\{V \in G(2,1) : \as P_V X = \ubd X/(1-\theta)\} \subseteq \{V \in G(2,1) : \as P_V X > u_\alpha\}$. Theorem \ref{thm:tube-lower-spectrum-alpha} gives $\hd \{V \in G(2,1) : \as P_V X = \ubd X/(1-\theta)\} \leq \alpha$ for every $\alpha \in (0,1)$, and letting $\alpha \downarrow 0$ proves the claim.  We emphasise that although the Hausdorff dimension of the exceptional set in \eqref{eq:zerodim} is zero, it is not necessarily empty, even under the assumption that $\as X =\ubd X < \ubd X/(1-\theta)$; see \cite[Theorem 3.4.12]{jon:book}.

\subsection{An application to self-similar sets}

For self-similar sets in the line, there is a clean dichotomy for Assouad dimension proved in \cite{fhor}: it either coincides with the Hausdorff dimension, or it is equal to 1. There is interest in determining the quasi-Assouad and Assouad spectrum of self-similar sets but this is still open \cite[Question 17.5.3]{jon:book}.  It is widely conjectured that in fact the quasi-Assouad dimension (and therefore Assouad spectrum) always  coincides with the Hausdorff dimension and this has been verified in some classes, including systems defined by algebraic parameters with no exact overlaps \cite[Theorem 7.3.1]{jon:book}.  By realising parametrized families of self-similar sets as projections of a common set, we are able to provide some evidence in the direction of this conjecture, that is, we give some non-trivial upper bounds on the Assouad spectrum of certain self-similar sets in a generic sense.

More precisely, projecting a planar self-similar set along a parametrized family of lines yields a one-parameter family of self-similar sets on $\R$. Let $I \subset \R$ be a bounded interval. For each $t \in I$, let $K_t$ be the self-similar set associated with the  homogeneous iterated function system
\begin{equation*}
  \Psi_t \ = \ \{\psi_{i,t}(x) = r x + a_i + t b_i\}_{i \in \{1,\ldots,N\}},
\end{equation*}
where $r \in (0,1)$ and $a_i,b_i \in \R$. Let $X \subset \R^2$ be the self-similar set  associated with the maps
\begin{equation*}
  F_i(x,y) \ = \ r(x,y) + (a_i,b_i)
\end{equation*}
for $i \in \{1,\ldots,N\}$, and $s = \log N/\log r^{-1}$ be the similarity dimension of $X$. Let $\mathcal{L}^d$ denote Lebesgue measure on $\R^d$.

\begin{thm} \label{thm:parametrised-spectrum}
  Suppose that, with the notation above, the planar self-similar set $X$ has quasi-Assouad dimension $\qad X \leq s$. Then, for every $\theta \in (0,1)$ and $\alpha \in [0,1)$,
  \begin{equation*}
    \hd \biggl\{ \ t \in I : \ \as K_t > \frac{1}{1+\alpha} \biggl( \frac{s}{1-\theta} + \alpha s \biggr)\biggr\} \ \leq \ \alpha.
  \end{equation*}
  Moreover,
  \begin{equation*}
    \as K_t \ \leq \ \frac{1}{2} \biggl( \frac{s}{1-\theta} + s \biggr)
  \end{equation*}
  for $\mathcal{L}^1$-almost all $t \in I$.
\end{thm}

See Section \ref{sec:parametrised-self-similar} for the proof. The hypothesis $\qad X \leq s$ is used only to simplify the bounds in Theorems \ref{thm:tube-lower-spectrum} and \ref{thm:tube-lower-spectrum-alpha}. The condition holds, for example, when $X$ satisfies the weak separation condition. However, taking the union of the exceptional sets for a sequence $\theta \uparrow 1$ gives no finite almost-sure bound for $\qad K_t$, because the upper threshold in Theorem \ref{thm:parametrised-spectrum} grows like $(1-\theta)^{-1}$ as $\theta \uparrow 1$. Sending $\alpha \downarrow 0$ gives
\begin{equation*}
  \hd \biggl\{ \ t \in I : \ \as K_t \ = \ \frac{s}{1-\theta} \biggr\} \ = \ 0
\end{equation*}
whenever $s = \log N/\log r^{-1} > 0$.

\begin{example} \label{ex:three-maps}
  Fix $0 < r < \frac12$, write $s = \log 3/\log r^{-1}$, and let $F_{\mathbf{b}} \subset \R$ be the homogeneous self-similar set obtained from three orientation-preserving similitudes
  \begin{equation*}
    f_i(x) \ = \ r x + b_i
  \end{equation*}
  with translations $\mathbf{b} = (b_1,b_2,b_3)$. It follows from Theorem \ref{thm:parametrised-spectrum} that, for every $\theta \in (0,1)$,
  \begin{equation} \label{eq:three-maps1}
    \as F_{\mathbf{b}} \ \leq \ \frac{1}{2} \biggl( \frac{s}{1-\theta} + s \biggr)
  \end{equation}
  for $\mathcal{L}^3$-almost all translations $\mathbf{b} = (b_1,b_2,b_3) \in \R^3$. Indeed, let $X \subset \R^2$ be the self-similar set generated by
  \begin{equation*}
    F_1(x,y) \ = \ r(x,y), \qquad F_2(x,y) \ = \ r(x,y) + (0,1-r), \qquad F_3(x,y) \ = \ r(x,y) + (1-r,0).
  \end{equation*}
  Since $X$ satisfies the strong separation condition, it is $s$-Ahlfors regular, so $\qad X = s$. Let $I = (0,1)$ and for each $t \in I$, let $K_t \subset \R$ be the homogeneous self-similar set generated by
  \begin{equation*}
    \psi_{1,t}(x) \ = \ r x, \qquad \psi_{2,t}(x) \ = \ r x + (1-r)t, \qquad \psi_{3,t}(x) \ = \ r x + (1-r).
  \end{equation*}
  By Theorem \ref{thm:parametrised-spectrum}, for each $\theta \in (0,1)$, we have
  \begin{equation*}
    \as K_t \ \leq \ \frac{1}{2} \biggl( \frac{s}{1-\theta} + s \biggr)
  \end{equation*}
  for $\mathcal{L}^1$-almost every $t \in I$. Let $E \subset I$ be the exceptional set where the above inequality fails.

  Suppose that $b_1 < b_2 < b_3$ and define
  \begin{equation*}
    a \ = \ b_1, \qquad d \ = \ b_3-b_1, \qquad t \ = \ \frac{b_2-b_1}{b_3-b_1} \ \in \ I.
  \end{equation*}
  If $h \colon \R \to \R$ is the homothety
  \begin{equation*}
    h(x) \ = \ \frac{d}{1-r}x + \frac{a}{1-r},
  \end{equation*}
  then
  \begin{equation*}
    f_1 \ = \ h \circ \psi_{1,t} \circ h^{-1}, \qquad f_2 \ = \ h \circ \psi_{2,t} \circ h^{-1}, \qquad f_3 \ = \ h \circ \psi_{3,t} \circ h^{-1}
  \end{equation*}
  and $F_{\mathbf{b}} = h(K_t)$. Since the Assouad spectrum is bi-Lipschitz invariant, translations $\mathbf{b} = (b_1,b_2,b_3)$ fail the bound \eqref{eq:three-maps1} if and only if $t \in E$.

  Define
  \begin{equation*}
    \Phi(a,d,t) \ = \ (a, a + d t, a + d).
  \end{equation*}
  Then $\Phi(\R \times (0,\infty) \times E)$ is exactly the exceptional set of ordered translation triples. Since $\L^1(E) = 0$, the set $\R \times (0,\infty) \times E$ has $\L^3$-measure zero. Moreover, $\Phi$ is locally Lipschitz on $\R \times (0,\infty) \times (0,1)$, so $\Phi(\R \times (0,\infty) \times E)$ is $\L^3$-null. For each permutation $\sigma$ of $\{1,2,3\}$, let $P_\sigma(b_1,b_2,b_3) = (b_{\sigma(1)}, b_{\sigma(2)}, b_{\sigma(3)})$. Then each $P_\sigma$ is an isometry of $\R^3$, so each set $P_\sigma(\Phi(\R \times (0,\infty) \times E))$ is $\L^3$-null. Since permuting the coordinates only relabels the similitudes, the full exceptional set is contained in the union of these finitely many sets and the set where two coordinates coincide, which is also $\L^3$-null. This proves \eqref{eq:three-maps1}.
\end{example}

\subsection{There is no Marstrand theorem for unbounded sets}

For unbounded sets, the Assouad spectrum of projections can attain arbitrary predetermined values on open sets of directions. The following result demonstrates that there is no Marstrand projection theorem for the Assouad spectrum at any fixed $\theta \in (0,1)$. A set $X \subseteq \R^d$ is \emph{uniformly discrete} if there exists $\delta > 0$ such that $|x-y| \geq \delta$ for all distinct $x,y \in X$. Every uniformly discrete subset of $\R^d$ is countable and closed, satisfies $\as X = 0$ for every $\theta \in (0,1)$, and hence has quasi-Assouad dimension zero.

\begin{thm} \label{main1}
  Let $\theta \in (0,1)$ and $s,t \in (0,1]$. Then there exists a uniformly discrete unbounded set $X \subseteq \mathbb{R}^2$ and non-empty disjoint open sets $I_s$ and $I_t$ in $G(2,1)$ such that
  \begin{enumerate}
    \item $\as P_V X = s$ for all $V \in I_s$, \vspace{2mm}
    \item $\as P_V X = t$ for all $V \in I_t$.
  \end{enumerate}
\end{thm}

See Section \ref{sec:projection-two-values} for the proof. The construction extends naturally to yield any finite or countable number of distinct values.

\begin{thm} \label{main1+}
  Let $\theta \in (0,1)$ and $s_j \in (0,1]$. Then there exists a uniformly discrete unbounded set $X \subseteq \mathbb{R}^2$ and non-empty disjoint open sets $\{I_j\}_{j\in \N}$  in $G(2,1)$ such that $\as P_V X = s_j$ for all $V \in I_j$.
\end{thm}

See Section \ref{sec:projection-infinite-values} for the proof. The preceding two theorems fix the parameter $\theta \in (0,1)$. The next result provides a single set that attains distinct values simultaneously for all $\theta \in (0,1)$.

\begin{thm} \label{maintheta}
  Let $0 < s < t \leq 1$. Then there exists a uniformly discrete unbounded set $X \subseteq \mathbb{R}^2$ and non-empty disjoint open sets $I_s$ and $I_t$ in $G(2,1)$ such that
  \begin{enumerate}
    \item $\as P_V X = s$ for all $\theta \in (0,1)$ and all $V \in I_s$, \vspace{2mm}
    \item $\as P_V X = t$ for all $\theta \in (0,1)$ and all $V \in I_t$.
  \end{enumerate}
\end{thm}

See Section \ref{sec:projection-simultaneous-theta} for the proof. Theorem \ref{maintheta} immediately implies that there is no Marstrand projection theorem for the quasi-Assouad dimension.

\begin{cor} \label{mainquasi}
  Let $0 < s < t \leq 1$. Then there exists a uniformly discrete unbounded set $X \subseteq \mathbb{R}^2$ and non-empty disjoint open sets $I_s$ and $I_t$ in $G(2,1)$ such that
  \begin{enumerate}
    \item $\qad P_V X = s$ for all $V \in I_s$, \vspace{2mm}
    \item $\qad P_V X = t$ for all $V \in I_t$.
  \end{enumerate}
\end{cor}

The upper box dimension of an unbounded set admits two natural extensions, and Marstrand's projection theorem holds for neither, although the obstructions are of a different nature. For bounded sets, $\ubd P_V X \leq \ubd X$ for every $V \in G(2,1)$, since orthogonal projections are $1$-Lipschitz. Under the usual extension $\ubd X = \sup_{n \in \N} \ubd(X \cap B(0,n))$, this may fail for unbounded sets: projection can strictly increase the upper box dimension. Indeed, the integer lattice $\mathbb{Z}^2$ is uniformly discrete, and each $\mathbb{Z}^2 \cap B(0,n)$ is finite, so $\ubd \mathbb{Z}^2 = 0$. On the other hand, for every line $V \in G(2,1)$ of irrational slope the image $P_V \mathbb{Z}^2$ is a dense subgroup of $V$, so $P_V \mathbb{Z}^2 \cap B(0,n)$ is dense in $B(0,n)$ and $\ubd(P_V \mathbb{Z}^2 \cap B(0,n)) = 1$ for every $n \in \N$. Hence $\ubd P_V \mathbb{Z}^2 = 1$ for $\gamma_{2,1}$-almost every $V \in G(2,1)$, even though $\ubd \mathbb{Z}^2 = 0$. Here $\ubd P_V \mathbb{Z}^2$ is still almost surely constant; the obstruction is that this constant exceeds $\ubd \mathbb{Z}^2$.

The usual extension nevertheless retains one Marstrand-type inequality. For $X \subseteq \mathbb{R}^d$ and an integer $1 \leq k < d$, extend the upper box dimension profile to $X$ by $\overline{\mbox{\rm dim}}_{\rm B}^k X = \sup_{n \in \N} \overline{\mbox{\rm dim}}_{\rm B}^k(X \cap B(0,n))$, in parallel with the usual extension of $\ubd$. Since $P_V(X \cap B(0,n)) \subseteq P_V X \cap B(0,n)$, monotonicity of the upper box dimension gives $\ubd(P_V X \cap B(0,n)) \geq \ubd P_V(X \cap B(0,n))$, and by Theorem \ref{mainA}(i) the right-hand side equals $\overline{\mbox{\rm dim}}_{\rm B}^k(X \cap B(0,n))$ for $\gamma_{d,k}$-almost every $V$. Intersecting these full-measure sets over $n \in \N$ and taking the supremum yields
\begin{equation*}
  \ubd P_V X \ \geq \ \overline{\mbox{\rm dim}}_{\rm B}^k X
\end{equation*}
for $\gamma_{d,k}$-almost every $V \in G(d,k)$. The integer lattice shows that the inequality can be strict, since $\overline{\mbox{\rm dim}}_{\rm B}^k \mathbb{Z}^2 = 0$ while $\ubd P_V \mathbb{Z}^2 = 1$ for almost every $V$.

The Wang--Li extension $\ubdzero X = \lim_{\theta \downarrow 0} \as X$ fails more sharply, and Theorem \ref{maintheta} yields the following analogue of Corollary \ref{mainquasi}.

\begin{cor} \label{mainwangli}
  Let $0 < s < t \leq 1$. Then there exists a uniformly discrete unbounded set $X \subseteq \mathbb{R}^2$ and non-empty disjoint open sets $I_s$ and $I_t$ in $G(2,1)$ such that
  \begin{enumerate}
    \item\label{it:wangli-s} $\ubdzero P_V X = s$ for all $V \in I_s$, \vspace{2mm}
    \item\label{it:wangli-t} $\ubdzero P_V X = t$ for all $V \in I_t$.
  \end{enumerate}
\end{cor}


In contrast to the lattice example, $\ubdzero P_V X$ is not almost surely constant: it takes the distinct values $s$ and $t$ on the open sets $I_s$ and $I_t$, each of which has positive $\gamma_{2,1}$-measure.

In view of Theorem \ref{thm:tube-lower-spectrum}, the situation for bounded planar sets is sharply different from the unbounded case treated in Theorems \ref{main1}, \ref{main1+}, and \ref{maintheta}, and it is natural to ask whether the Assouad spectrum of projections of bounded planar sets satisfies a Marstrand-type identity.

\begin{ques} \label{ques:bounded-marstrand-spectrum}
  Is there a Marstrand projection theorem for the Assouad spectrum of bounded planar sets?
\end{ques}

The tube-counting bound of Theorem \ref{thm:tube-lower-spectrum} is far from sharp. The matching lower bound from Theorem \ref{main2} goes through a different invariant, the dimension profile $\dim_{\textup{A}}^{\theta,k} X$, and the two bounds need not meet. A genuine Marstrand-type identity---almost-sure equality between $\as P_V X$ and some intrinsic spectrum-profile of $X$---is consistent with everything currently proved.

Throughout the paper, for $A, B > 0$, we use the notation $A \lesssim B$ to indicate that there exists an absolute constant $C > 0$ such that $A \leq C B$. We write $A \approx B$ to denote the two-sided inequality $A \lesssim B \lesssim A$.

\section{Almost sure lower bound for the Assouad spectrum:~proof of Theorem \ref{main2}} \label{sec:almost-sure-lower-bound}

In this section we prove Theorem \ref{main2} and derive Corollary \ref{cor:assouad-spectrum-profile-bound}. The proof of the theorem turns the $k$-dimensional Assouad spectrum profile into an almost-sure lower bound for projections; the proof of the corollary then compares the $d$- and $k$-dimensional capacity estimates to obtain the explicit lower bound.

To prove Theorem \ref{main2}, we show that every $\alpha < \dim_{\textup{A}}^{\theta,k} X$ is also a lower bound for $\as P_V X$ for $\gamma_{d,k}$-almost all $V \in G(d,k)$. The proof relies on energy estimates and Markov's inequality to control the covering numbers of the projected sets. Let $\alpha < \dim_{\textup{A}}^{\theta,k} X$. By the definition of the $k$-dimensional Assouad spectrum profile, for any $C > 0$, there exist $x \in X$ and $0 < R < 1$ such that
\begin{equation} \label{star}
  C_{R^{1/\theta}}^k(X \cap B(x, R)) \ \geq \ C \biggl(\frac{R}{R^{1/\theta}}\biggr)^\alpha.
\end{equation}
In the proof, we will later select a sequence of constants $C > 0$, but for now, we fix $C$ and let $x \in X$ and $0 < R < 1$ be such that \eqref{star} holds. Set $E = \overline{X \cap B(x, R)}$. Since the capacity of a bounded set is defined as the capacity of its closure, \eqref{star} yields
\begin{equation*}
  C_{R^{1/\theta}}^k(E) \ \geq \ C \biggl(\frac{R}{R^{1/\theta}}\biggr)^\alpha.
\end{equation*}
Let $\mu \in \mathcal{M}(E)$ be an equilibrium measure associated with the capacity $C_{R^{1/\theta}}^k(E)$. It follows that
\begin{equation*}
  \iint \phi_{R^{1/\theta}}^k(w- y) \,\mathrm{d}\mu(w) \,\mathrm{d}\mu(y) \ = \ C_{R^{1/\theta}}^k(E)^{-1} \ \leq \ C^{-1} \biggl(\frac{R}{R^{1/\theta}}\biggr)^{-\alpha}.
\end{equation*}
Let $\mu_V = (P_V)_* \mu \in \mathcal{M}(P_V E)$ be the pushforward measure under the projection $P_V$. Then
\begin{align*}
  \int (\mu_V \times \mu_V) &(\{(z, u) \in B( P_V  x, R)^2 : |z - u| \leq R^{1/\theta} \})\,\mathrm{d}\gamma_{d,k}(V) \\
  &= \ \int (\mu \times \mu)(\{(w, y) \in B(x, R)^2 : | P_V w -  P_V  y| \leq R^{1/\theta} \})\,\mathrm{d}\gamma_{d,k}(V) \\
  &= \ \iint \gamma_{d,k}(\{V \in G(d, k) : | P_V w -  P_V  y| \leq R^{1/\theta}\})\,\mathrm{d}\mu(w)\,\mathrm{d}\mu(y) \\
  &\leq \ C_{d,k} \iint \phi_{R^{1/\theta}}^k(w- y) \,\mathrm{d}\mu(w)\,\mathrm{d}\mu(y) \\
  &\leq \ C^{-1} C_{d,k} \biggl(\frac{R}{R^{1/\theta}}\biggr)^{-\alpha},
\end{align*}
where $C_{d,k}$ is a constant depending only on $d$ and $k$. By Markov's inequality,
\begin{align*}
  \gamma_{d,k}\biggl( &\biggl\{V : (\mu_V \times \mu_V)(\{(z, u) \in B( P_V  x, R)^2 : |z - u| \leq R^{1/\theta}\}) > C^{-1/2} C_{d,k} \biggl(\frac{R}{R^{1/\theta}}\biggr)^{-\alpha} \biggr\}\biggr) \\
  &\leq \ C^{1/2} C_{d,k}^{-1} \biggl(\frac{R}{R^{1/\theta}}\biggr)^{\alpha} \int (\mu_V \times \mu_V)(\{(z, u) \in B( P_V  x, R)^2 : |z - u| \leq R^{1/\theta} \})\,\mathrm{d}\gamma_{d,k}(V) \\
  &\leq \ C^{-1/2}
\end{align*}
Thus, by \cite[Lemma 2.2]{falconerprofile}, there exists a constant $c_d$ depending only on $d$ such that
\begin{equation} \label{eq:lower-bound-gamma-measure}
  \gamma_{d,k}\biggl(\biggl\{V : N_{R^{1/\theta}}(P_V X \cap B( P_V x , R)) < C^{1/2} c_d C_{d,k}^{-1} \biggl(\frac{R}{R^{1/\theta}}\biggr)^\alpha\biggr\}\biggr) \ \leq \ C^{-1/2}.
\end{equation}
Indeed, since $P_V E$ is compact, applying \cite[Lemma 2.2]{falconerprofile} to $\mu_V$ yields a lower bound for $N_{R^{1/\theta}}(P_V E)$. Also, $P_V(X \cap B(x,R)) \subseteq P_V X \cap B(P_V x,R)$ and $P_V(X \cap B(x,R)) \subseteq P_V E \subseteq \overline{P_V(X \cap B(x,R))}$, so
\begin{equation*}
  N_r(P_V E) \ = \ N_r(P_V(X \cap B(x,R))) \ \leq \ N_r(P_V X \cap B(P_V x,R))
\end{equation*}
for every $r > 0$. Hence this lower bound carries over to $N_{R^{1/\theta}}(P_V X \cap B(P_V x,R))$, which gives \eqref{eq:lower-bound-gamma-measure}.

We now select the sequence of constants $C > 0$. For each $i \in \mathbb{N}$, set $C_i = 2^{2i}$ and choose corresponding $x_i \in X$ and $0 < R_i < 1$ satisfying \eqref{star}. Define
\begin{equation*}
  G_i \ = \ \biggl\{V : N_{R_i^{1/\theta}}(P_V X \cap B( P_V x_i, R_i)) < 2^i c_d C_{d,k}^{-1} \biggl(\frac{R_i}{R_i^{1/\theta}}\biggr)^\alpha\biggr\}
\end{equation*}
and
\begin{equation*}
  G \ = \ \bigcap_{l = 1}^\infty \bigcup_{i = l}^\infty G_i.
\end{equation*}
Since $\sum_{i=1}^\infty \gamma_{d,k}(G_i) \leq \sum_{i=1}^\infty 2^{-i} < \infty$ by \eqref{eq:lower-bound-gamma-measure}, the Borel-Cantelli lemma implies $\gamma_{d,k}(G) = 0$. If $V \in G(d,k) \setminus G$, then there exists $i_0 \in \mathbb{N}$ such that $V \notin G_i$ for all $i \geq i_0$. Hence, we have $P_V x_i \in P_V X$ and
\begin{equation*}
  N_{R_i^{1/\theta}}(B( P_V x_i, R_i) \cap P_V X) \ \geq \ 2^i c_d C_{d,k}^{-1} \biggl(\frac{R_i}{R_i^{1/\theta}}\biggr)^\alpha
\end{equation*}
for all $i \geq i_0$. Thus no constant works in the definition of $\as(P_V X)$ at the exponent $\alpha$, and so $\as(P_V X) \geq \alpha$ for all $V \in G(d,k) \setminus G$. If $\dim_{\textup{A}}^{\theta,k} X = 0$, the theorem is immediate since $\as(P_V X) \geq 0$ for all $V \in G(d,k)$. Otherwise, applying the argument above to each rational $\alpha \in \Q \cap (0,\dim_{\textup{A}}^{\theta,k} X)$ and intersecting the corresponding full-measure sets yields $\as(P_V X) \geq \dim_{\textup{A}}^{\theta,k} X$ for $\gamma_{d,k}$-almost all $V \in G(d,k)$, as required.

We now derive Corollary \ref{cor:assouad-spectrum-profile-bound} from Theorem \ref{main2} by proving the corresponding lower bound for the $k$-dimensional Assouad spectrum profile at the parameter $\theta'$. If $\as X = 0$, the assertion is immediate from Theorem \ref{main2}, so assume that $\as X > 0$ and fix $\alpha \in (0,\as X)$. Let $C > 1$. By the capacity formulation \eqref{eq:assouad-spectrum-characterization} of the Assouad spectrum, derived from \cite[Corollary 2.4]{falconerprofile}, there are $x \in X$ and $0 < R < 1$ such that, with $E = \overline{X \cap B(x,R)}$ and $\rho = R^{1/\theta}$,
\begin{equation*}
  C_\rho^d(E) \ \geq \ C \biggl(\frac{R}{\rho}\biggr)^\alpha.
\end{equation*}
Let $\mu \in \mathcal{M}(E)$ be an equilibrium measure associated with $C_\rho^d(E)$. Then
\begin{equation*}
  \iint \phi_\rho^d(w-y)\,\mathrm{d}\mu(w)\,\mathrm{d}\mu(y) \ \leq \ C^{-1} R^{(1/\theta-1)\alpha}.
\end{equation*}
Put $r = R^{1/\theta'}$, so that $r \leq \rho$. For fixed $w \in E$, split the integral below into the ranges $|w-y| < \rho$ and $|w-y| \geq \rho$. On the first range, $\phi_r^k(w-y) \leq \phi_\rho^d(w-y)$, while on the second range $\phi_r^k(w-y) = (r/\rho)^k\phi_\rho^d(w-y)^{k/d}$. Since $\mu$ is a probability measure and $k/d < 1$, Jensen's inequality applied to the second range gives
\begin{equation*}
  \int \phi_r^k(w-y)\,\mathrm{d}\mu(y) \ \leq \ \int \phi_\rho^d(w-y)\,\mathrm{d}\mu(y) + \biggl(\frac{r}{\rho}\biggr)^k \biggl(\int \phi_\rho^d(w-y)\,\mathrm{d}\mu(y)\biggr)^{k/d}.
\end{equation*}
Integrating in $w$ and using Jensen's inequality once more gives
\begin{align*}
  \iint \phi_r^k(w-y)\,\mathrm{d}\mu(w)\,\mathrm{d}\mu(y) \ &\leq \ C^{-1} R^{(1/\theta-1)\alpha} + C^{-k/d} R^{k(1/\theta'-1/\theta)+(1/\theta-1)\alpha k/d} \\
  &\leq \ 2C^{-k/d} R^{(1/\theta-1)\alpha}.
\end{align*}
In the last inequality we used the definition of $\theta'$ and the inequality $\alpha < \as X$, which imply
\begin{equation*}
  k\biggl(\frac{1}{\theta'}-\frac{1}{\theta}\biggr)+\biggl(\frac{1}{\theta}-1\biggr)\frac{\alpha k}{d} \ \geq \ \biggl(\frac{1}{\theta}-1\biggr)\alpha.
\end{equation*}
Since
\begin{equation*}
  R^{-(1/\theta-1)\alpha} \ = \ \biggl(\frac{R}{R^{1/\theta'}}\biggr)^{\alpha/(1+(1/k-1/d)\as X)},
\end{equation*}
it follows that
\begin{equation*}
  C_{R^{1/\theta'}}^k(E) \ \geq \ \frac{1}{2} C^{k/d} \biggl(\frac{R}{R^{1/\theta'}}\biggr)^{\alpha/(1+(1/k-1/d)\as X)}.
\end{equation*}
Since $C>1$ was arbitrary, the definition of the $k$-dimensional Assouad spectrum profile gives
\begin{equation*}
  \dim_{\textup{A}}^{\theta',k} X \ \geq \ \frac{\alpha}{1+(\frac{1}{k}-\frac{1}{d})\as X}.
\end{equation*}
Letting $\alpha \uparrow \as X$ and applying Theorem \ref{main2} with $\theta'$ in place of $\theta$ proves the claim.

\section{Tube-counting obstructions for spectra of projections:~proof of Theorem \ref{thm:tube-lower-spectrum}} \label{sec:tube-counting-obstructions}

To prove Theorem \ref{thm:tube-lower-spectrum}, we extract a well-separated subset of $X$ by carefully selecting points from tubes corresponding to exceptional projection directions, thereby bounding $\as X$ from below. Replacing $X$ by its closure $\overline{X}$ does not change $\as X$, $\ubd X$, or $\as P_V X$; therefore, we may assume that $X$ is compact. If
\begin{equation*}
  \esssup_{V \in G(2,1)} \as P_V X \ \leq \ \as X,
\end{equation*}
then, since $\as X \leq \ubd X/(1-\theta)$ for bounded sets, we have
\begin{equation*}
  \esssup_{V \in G(2,1)} \as P_V X \ \leq \ \frac{1}{2} \biggl( \frac{\ubd X}{1-\theta} + \as X \biggr),
\end{equation*}
and the conclusion follows immediately. It therefore remains to consider the case
\begin{equation*}
  \as X \ < \ \esssup_{V \in G(2,1)} \as P_V X.
\end{equation*}
In this case there exist $s,t > 0$ satisfying
\begin{equation*}
  \as X \ < \ s \ < \ t \ < \ \esssup_{V \in G(2,1)} \as P_V X.
\end{equation*}
Let $S_k = 2^{-k}$ for all $k \in \N$. Fix a constant $C_\theta \geq 1$ such that every interval of radius $(2R)^{1/\theta}$ can be covered by at most $C_\theta$ intervals of radius $R^{1/\theta}$ for all $0 < R < 1$. For each $k,n \in \N$, define
\begin{align*}
  J_{k,n} \ = \ \biggl\{V \in G(2,1) &: \text{   there exists } x \in P_V X \text{ such that} \\
  &N_{S_k^{1/\theta}} \big( P_V X \cap B(x,S_k) \big) \ > \ n C_\theta^{-1} \biggl(\frac{S_k}{S_k^{1/\theta}}\biggr)^t \biggr\}.
\end{align*}
For fixed $k,n$, set
\begin{equation*}
  m_{k,n} \ = \ \biggl\lfloor n C_\theta^{-1} \biggl(\frac{S_k}{S_k^{1/\theta}}\biggr)^t \biggr\rfloor + 1.
\end{equation*}
Since $P_V X \cap B(x,S_k)$ is a compact subset of the line $V$, the condition in the definition of $J_{k,n}$ is equivalent to the existence of $m_{k,n}$ points in $P_V X \cap B(x,S_k)$ that are pairwise more than $2S_k^{1/\theta}$ apart. Indeed, for a compact set $K$ in the line $V$ and $r > 0$, the inequality $N_r(K) > A$ is equivalent to the existence of $\lfloor A \rfloor + 1$ points of $K$ that are pairwise more than $2r$ apart. The reverse implication is immediate, and the forward implication follows from the standard interval covering argument on the line. Hence $J_{k,n}$ is the projection onto $G(2,1)$ of the Borel set of tuples
\begin{equation*}
  (V,z,y_1,\ldots,y_{m_{k,n}}) \ \in \ G(2,1) \times X^{m_{k,n}+1}
\end{equation*}
satisfying
\begin{equation*}
  |P_V y_i - P_V z| \ \leq \ S_k
\end{equation*}
for all $i \in \{1,\ldots,m_{k,n}\}$ and
\begin{equation*}
  |P_V y_i - P_V y_j| \ > \ 2S_k^{1/\theta}
\end{equation*}
for all distinct $i,j \in \{1,\ldots,m_{k,n}\}$. Therefore each set $J_{k,n}$ is analytic, and hence $\gamma_{2,1}$-measurable. Set
\begin{equation*}
  I \ = \ \bigcap_{n = 1}^\infty \limsup_{k \to \infty} J_{k,n}.
\end{equation*}
Then $I$ is measurable. For each $m,k,n \in \N$, let
\begin{align*}
  J_{k,n}^{(m)} \ = \ \biggl\{V \in G(2,1) &: \text{   there exists } x \in P_V X \text{ such that} \\
  &N_{S_k^{1/\theta}} \big( P_V X \cap B(x,S_k) \big) \ > \ n C_\theta^{-1} \biggl(\frac{S_k}{S_k^{1/\theta}}\biggr)^{t+1/m} \biggr\},
\end{align*}
and set
\begin{equation*}
  E_m \ = \ \bigcap_{n = 1}^\infty \limsup_{k \to \infty} J_{k,n}^{(m)}.
\end{equation*}
The same projection argument, with $t + 1/m$ in place of $t$, shows that each set $J_{k,n}^{(m)}$ is analytic, and hence each set $E_m$ is measurable. If $\as P_V X > t$, then there exists $m \in \N$ such that $\as P_V X > t + 1/m$. In this case, no constant works in the definition of the Assouad spectrum at the exponent $t + 1/m$. Thus, for any $n \in \N$, there are arbitrarily small scales $R > 0$ and points $x_R \in P_V X$ for which the covering number $N_{R^{1/\theta}}(P_V X \cap B(x_R, R))$ exceeds $n(R/R^{1/\theta})^{t+1/m}$. Choosing $k \in \N$ such that $S_{k+1} \le R < S_k$, we observe that an $S_k^{1/\theta}$-interval can be covered by at most $C_\theta$ intervals of radius $R^{1/\theta}$, yielding $N_{S_k^{1/\theta}}(P_V X \cap B(x_R, S_k)) > n C_\theta^{-1}(S_k/S_k^{1/\theta})^{t+1/m}$. Since this holds for arbitrarily small $R$, we have $V \in \limsup_{k \to \infty} J_{k,n}^{(m)}$ for all $n \in \N$, and hence $V \in E_m$. Conversely, if $V \in E_m$ for some $m \in \N$ and there were a constant $C > 0$ such that
\begin{equation*}
  N_{R^{1/\theta}} \big( P_V X \cap B(x,R) \big) \ \leq \ C \biggl(\frac{R}{R^{1/\theta}}\biggr)^{t+1/m}
\end{equation*}
for all $x \in P_V X$ and all $0 < R < 1$, then choosing $n \in \N$ with $n > C C_\theta$ would contradict the fact that $V \in J_{k,n}^{(m)}$ for infinitely many $k \in \N$. Hence $\as P_V X > t$. Therefore
\begin{equation*}
  \{V \in G(2,1) : \as P_V X > t\} \ = \ \bigcup_{m = 1}^\infty E_m.
\end{equation*}
Since $t < \esssup_{V \in G(2,1)} \as P_V X$, the set on the left-hand side has positive $\gamma_{2,1}$-measure, and therefore some $E_m$ has positive $\gamma_{2,1}$-measure. If $V \in E_m$, then for every $n \in \N$ there are infinitely many $k \in \N$ such that $V \in J_{k,n}^{(m)}$. Since $0 < S_k < 1$ and $\theta \in (0,1)$, we have $S_k / S_k^{1/\theta} > 1$, and therefore $J_{k,n}^{(m)} \subseteq J_{k,n}$ because $t + 1/m > t$. Thus every direction in $E_m$ belongs to $I$, and so $\gamma_{2,1}(I) > 0$. The parameter $n$ in the definition of $J_{k,n}$ serves to neutralize the arbitrary constant $C$ in the definition of the Assouad spectrum. For the measure estimates in the remainder of the proof, setting $n=1$ is sufficient. For brevity, write $J_k = J_{k,1}$.

By the definition of $I$, if $V \in I$, then $V \in \limsup_{k \to \infty} J_{k,n}$ for every $n \in \N$. On the other hand, for every $n \in \N$ we have $J_{k,n} \subseteq J_k$. This follows since the condition
\begin{equation*}
  N_{S_k^{1/\theta}}( P_V X \cap B(x,S_k) ) > n C_\theta^{-1} \biggl(\frac{S_k}{S_k^{1/\theta}}\biggr)^t
\end{equation*}
defining $J_{k,n}$ becomes strictly weaker as $n$ decreases. Hence $V \in \limsup_{k \to \infty} J_k$ and, consequently,
\begin{equation*}
  I \ \subseteq \ \limsup_{k \to \infty} J_k.
\end{equation*}
By the Borel-Cantelli lemma, this implies
\begin{equation*}
  \sum_{k = 1}^\infty \gamma_{2,1}(J_k) \ = \ \infty,
\end{equation*}
and therefore, for each $\eta > 0$, there are infinitely many $k \in \N$ such that
\begin{equation*}
  \gamma_{2,1}(J_k) \ \geq \ S_k^\eta.
\end{equation*}
Fix a sufficiently large such $k$ and set $R = S_k$. Since there are infinitely many such $k$ and $R \downarrow 0$, we may assume that $2R < 1$. Choose
\begin{equation} \label{eq:choice-of-beta}
  0 \ < \ \eta \ < \ \beta \ < \ \frac{(t-s)(1-\theta)}{\theta},
\end{equation}
and consider a maximal $R^\beta$-separated subset $\mathcal{E} \subseteq J_k$. Identifying $G(2,1)$ with a circle equipped with normalised arc-length metric, every ball of radius $R^\beta$ has $\gamma_{2,1}$-measure comparable to $R^\beta$. Hence
\begin{equation*}
  \# \mathcal{E} \ \gtrsim \ \gamma_{2,1}(J_k) R^{-\beta} \ \geq \ R^{-(\beta-\eta)},
\end{equation*}
and also $\# \mathcal{E} \lesssim R^{-\beta}$. For each $V \in \mathcal{E}$ there is $x_V \in P_V X$ such that
\begin{equation} \label{eq:uniform-lower-spectrum-constant}
  N_{R^{1/\theta}} \big( P_V X \cap B(x_V,R) \big) \ > \ C_\theta^{-1} \biggl(\frac{R}{R^{1/\theta}}\biggr)^t.
\end{equation}
Set $r = (2R)^{1/\theta}$. Since $\theta \in (0,1)$ and $R$ is sufficiently small, we also have $2r < R$. For each $V \in \mathcal{E}$, since each ball of radius $r$ in the line $V$ can be covered by at most $C_\theta$ balls of radius $R^{1/\theta}$, \eqref{eq:uniform-lower-spectrum-constant} yields
\begin{equation*}
  N_r \big( P_V X \cap B(x_V,R) \big) \ > \ C_\theta^{-2} \biggl(\frac{R}{R^{1/\theta}}\biggr)^t.
\end{equation*}
Let $T_V = P_V^{-1}B(x_V,R)$ be the associated $R$-tube in the plane, and choose a maximal subset $F_V \subseteq P_V X \cap B(x_V,R)$ whose points are pairwise more than $2r$ apart. Since the $2r$-balls centred at points of $F_V$ cover $P_V X \cap B(x_V,R)$ and each such ball in the line $V$ can be covered by two balls of radius $r$, we have
\begin{equation*}
  \# F_V \ \geq \ \frac{1}{2} N_r \big( P_V X \cap B(x_V,R) \big) \ > \ \frac{1}{2} C_\theta^{-2} \biggl(\frac{R}{R^{1/\theta}}\biggr)^t.
\end{equation*}
For each $y \in F_V$, choose a point $z_y \in X$ with $P_V z_y = y$, and set $E_V = \{z_y : y \in F_V\}$. Since $P_V$ is $1$-Lipschitz, distinct points of $E_V$ are more than $2r$ apart. Moreover,
\begin{equation*}
  \# E_V \ > \ \frac{1}{2} C_\theta^{-2} \biggl(\frac{R}{R^{1/\theta}}\biggr)^t.
\end{equation*}
By passing to a subset of $\mathcal{E}$ of cardinality at least $\frac{1}{2} \# \mathcal{E}$ and relabeling, we may assume that all directions in $\mathcal{E}$ are contained in an angular interval of length $\pi/2$. Since the bounds above for $\# \mathcal{E}$ are unchanged up to constants, we still have $\# \mathcal{E} \gtrsim R^{-(\beta-\eta)}$ and $\# \mathcal{E} \lesssim R^{-\beta}$. Order the directions in $\mathcal{E}$ by the angle that they make with the $x$-axis, and write $\mathcal{E} = \{V_1,\ldots,V_M\}$ with $M = \# \mathcal{E}$. Since $\mathcal{E}$ is $R^\beta$-separated and there is no cyclic wrap-around inside this interval, the angle between $V_m$ and $V_{m-i}$ is at least $i R^\beta$ for every $1 \leq i < m \leq M$. For each $V \in \mathcal{E}$, let $T_V[2r]$ be the $2r$-neighbourhood of $T_V$. Because $2r < R$, enlarging one tube by $2r$ changes its width by at most a constant factor. For each $j \in \{1,\ldots,M\}$, choose a unit vector $e_j \in S^1 \cap V_j$ and a scalar $a_j \in \R$ such that
\begin{equation*}
  T_{V_j} \ = \ \{z \in \R^2 : |z \cdot e_j - a_j| \leq R\}.
\end{equation*}
Let $\alpha_{i,m}$ denote the angle between $V_m$ and $V_{m-i}$. Then $\alpha_{i,m} \geq i R^\beta$. Since $T_{V_{m-i}}[2r]$ also has width comparable to $R$, the set $T_{V_m} \cap T_{V_{m-i}}[2r]$ is contained in a parallelogram whose side lengths are at most a constant multiple of $R$ and $R/\sin \alpha_{i,m}$. Indeed, both strips have width comparable to $R$, and since they meet at angle $\alpha_{i,m}$, their intersection is contained in such a parallelogram. A ball of radius $R$ covers a segment of length comparable to $R$ in the long direction, so the parallelogram can be covered by at most $C (\sin \alpha_{i,m})^{-1}$ balls of radius $R$. As $\alpha_{i,m} \in [0,\pi/2]$, we have $\sin \alpha_{i,m} \gtrsim \alpha_{i,m} \gtrsim i R^\beta$, and so there exists a constant $C_2 > 0$ such that, for every $1 \leq i < m \leq M$, the intersection $T_{V_m} \cap T_{V_{m-i}}[2r]$ can be covered by at most $C_2 i^{-1} R^{-\beta}$ balls of radius $R$. Since $s > \as X$, there exists a constant $C_s > 0$ such that
\begin{equation*}
  N_{\rho^{1/\theta}} \big( X \cap B(x,\rho) \big) \ \leq \ C_s \biggl(\frac{\rho}{\rho^{1/\theta}}\biggr)^s
\end{equation*}
for all $x \in X$ and all $0 < \rho < 1$. If $B(z,R)$ is one of the balls in the above cover and $B(z,R) \cap X \neq \emptyset$, choose $x \in B(z,R) \cap X$. Then $B(z,R) \cap X \subseteq B(x,2R) \cap X$, and so
\begin{equation*}
  N_r \big( B(z,R) \cap X \big) \ \leq \ N_r \big( B(x,2R) \cap X \big) \ \leq \ C_3 \biggl(\frac{R}{R^{1/\theta}}\biggr)^s,
\end{equation*}
where $C_3 > 0$ depends only on $X$, $s$, and $\theta$.

With the geometric estimates in place, we extract the required separated subset of $X$ by choosing points from the sets $E_{V_m}$ one tube at a time. When treating the $m$th tube, discard every point of $E_{V_m}$ that lies within distance $2r$ of a point selected from one of the previous tubes. Any discarded point must lie in $T_{V_m} \cap T_{V_j}[2r]$ for some $j < m$. Since distinct points of $E_{V_m}$ are more than $2r$ apart, each $r$-ball contains at most one point of $E_{V_m}$. Therefore the number of discarded points from the $m$th tube is at most
\begin{equation*}
  \sum_{i = 1}^{m-1} C_2 C_3 i^{-1} R^{-\beta} \biggl(\frac{R}{R^{1/\theta}}\biggr)^s \ \leq \ C_4 (\log m) R^{-\beta} \biggl(\frac{R}{R^{1/\theta}}\biggr)^s,
\end{equation*}
where $C_4 > 0$ is independent of $m$ and $R$. Hence the $m$th tube contributes at least
\begin{equation*}
  \frac{1}{2} C_\theta^{-2} \biggl(\frac{R}{R^{1/\theta}}\biggr)^t - C_4 (\log m) R^{-\beta} \biggl(\frac{R}{R^{1/\theta}}\biggr)^s
\end{equation*}
new points. Since $m \leq M \lesssim R^{-\beta}$, we have $\log m \lesssim |\log R|$. Moreover,
\begin{equation*}
  |\log R| R^{-\beta + (s-t)(1-1/\theta)} \ \to \ 0
\end{equation*}
as $R \downarrow 0$, because
\begin{equation*}
  -\beta + (s-t)(1-1/\theta) \ = \ -\beta + \frac{(t-s)(1-\theta)}{\theta} \ > \ 0
\end{equation*}
by \eqref{eq:choice-of-beta}. Hence, for all sufficiently small such $R$, the $m$th tube contributes at least
\begin{equation*}
  \frac{1}{4} C_\theta^{-2} \biggl(\frac{R}{R^{1/\theta}}\biggr)^t
\end{equation*}
new points. We therefore obtain a $2r$-separated subset of $X$ with cardinality at least
\begin{equation*}
  \gtrsim \ \sum_{m = 1}^{M} \biggl(\frac{R}{R^{1/\theta}}\biggr)^t \ \gtrsim \ R^{-(\beta-\eta)} \biggl(\frac{R}{R^{1/\theta}}\biggr)^t \ = \ (R^{-1/\theta})^{(1-\theta)t + (\beta-\eta)\theta}.
\end{equation*}
Since each ball of radius $r$ contains at most one point of this set, the same lower bound holds for $N_r(X)$. As $r = (2R)^{1/\theta}$ is comparable to $R^{1/\theta}$ and there are infinitely many such scales $R \downarrow 0$, we deduce that
\begin{equation*}
  \ubd X \ \geq \ (1-\theta)t + (\beta-\eta)\theta.
\end{equation*}
This lower bound is increasing in $\beta$, so letting $\beta \uparrow (t-s)(1-\theta)/\theta$ yields
\begin{equation*}
  t \ \leq \ \frac{1}{2} \biggl( \frac{\ubd X}{1-\theta} + s + \frac{\eta\theta}{1-\theta} \biggr).
\end{equation*}
Since this holds for every $\eta > 0$ with $\eta < (t-s)(1-\theta)/\theta$, every $t < \esssup_{V \in G(2,1)} \as P_V X$, and every $s > \as X$, we conclude that
\begin{equation*}
  \esssup_{V \in G(2,1)} \as P_V X \ \leq \ \frac{1}{2} \biggl( \frac{\ubd X}{1-\theta} + \as X \biggr)
\end{equation*}
as claimed.

\section{Exceptional set estimates for  spectra of projections:~proof of Theorem \ref{thm:tube-lower-spectrum-alpha}} \label{sec:tube-counting-obstructions-alpha}

The proof of Theorem \ref{thm:tube-lower-spectrum-alpha} is a straightforward modification of the argument used in Section \ref{sec:tube-counting-obstructions}. Let $X \subset \R^2$ be bounded and $\theta \in (0,1)$. Replacing $X$ by its closure $\overline{X}$ does not change $\as X$, $\ubd X$, or $\as P_V X$; therefore, we may assume that $X$ is compact. Let $u \in (\as X,\ubd X/(1-\theta)]$. We first establish that
\begin{equation} \label{eq:tube-superlevel-bound}
  \hd \{ V \in G(2,1) : \as P_V X > u\} \ \leq \ \frac{\ubd X/(1-\theta) - u}{u - \as X}.
\end{equation}
If $u = \ubd X/(1-\theta)$, then the exceptional set is empty by the trivial bound $\as P_V X \leq \ubd X/(1-\theta)$ for all $V \in G(2,1)$. So assume that $u < \ubd X/(1-\theta)$, and write
\begin{equation*}
  E_u \ = \ \{V \in G(2,1) : \as P_V X > u\}.
\end{equation*}
Assume for contradiction that
\begin{equation*}
  \frac{\ubd X/(1-\theta) - u}{u - \as X} \ < \ \hd E_u
\end{equation*}
and choose $\delta$ strictly between the two values. Since
\begin{equation*}
  \frac{\ubd X/(1-\theta) - u}{u - s} \ \to \ \frac{\ubd X/(1-\theta) - u}{u - \as X}
\end{equation*}
as $s \downarrow \as X$, we may choose $s \in (\as X,u)$ such that
\begin{equation} \label{eq:delta-choice}
  \delta \ > \ \frac{\ubd X/(1-\theta) - u}{u - s}.
\end{equation}
Let $S_k = 2^{-k}$ for all $k \in \N$, and fix a constant $C_\theta \geq 1$ such that every interval of radius $(2R)^{1/\theta}$ can be covered by at most $C_\theta$ intervals of radius $R^{1/\theta}$ for all $0 < R < 1$. For $k,n \in \N$, let
\begin{align*}
  J_{k,n}(u) \ = \ \biggl\{ \ V \in G(2,1) &: \text{   there exists } x \in P_V X \text{ such that} \\
  &N_{S_k^{1/\theta}} \big( P_V X \cap B(x,S_k) \big) \ > \ n C_\theta^{-1} \biggl(\frac{S_k}{S_k^{1/\theta}}\biggr)^u \biggr\},
\end{align*}
and let
\begin{equation*}
  I(u) \ = \ \bigcap_{n = 1}^\infty \limsup_{k \to \infty} J_{k,n}(u).
\end{equation*}
For $m,k,n \in \N$, let
\begin{align*}
  J_{k,n}^{(m)}(u) \ = \ \biggl\{ \ V \in G(2,1) &: \text{   there exists } x \in P_V X \text{ such that} \\
  &N_{S_k^{1/\theta}} \big( P_V X \cap B(x,S_k) \big) \ > \ n C_\theta^{-1} \biggl(\frac{S_k}{S_k^{1/\theta}}\biggr)^{u+1/m} \biggr\},
\end{align*}
and set
\begin{equation*}
  E_m(u) \ = \ \bigcap_{n = 1}^\infty \limsup_{k \to \infty} J_{k,n}^{(m)}(u).
\end{equation*}
Exactly as in the proof of Theorem \ref{thm:tube-lower-spectrum}, the sets $J_{k,n}(u)$ and $J_{k,n}^{(m)}(u)$ are analytic. Since countable unions and intersections of analytic sets are analytic, the sets $I(u)$ and $E_m(u)$ are analytic. The same argument as before also yields
\begin{equation*}
  E_u \ = \ \bigcup_{m = 1}^\infty E_m(u).
\end{equation*}
Hence there exists $m \in \N$ with $\hd E_m(u) > \delta$. By Frostman's lemma, there is a Radon probability measure $\sigma$ with $\spt \sigma \subset E_m(u)$ and
\begin{equation*}
  \sigma(B(V,\rho)) \ \lesssim \ \rho^\delta
\end{equation*}
for all $V \in G(2,1)$ and $\rho > 0$. Since analytic subsets of the compact metric space $G(2,1)$ are universally measurable, the same holds for $I(u)$ and $E_m(u)$. Since $\sigma$ is supported on $E_m(u) \subset I(u)$ and, just as before, $I(u) \subset \limsup_{k \to \infty} J_k(u)$ with $J_k(u) = J_{k,1}(u)$, we have $\sigma(\limsup_{k \to \infty} J_k(u)) = 1$. The Borel-Cantelli lemma therefore implies
\begin{equation*}
  \sum_{k = 1}^\infty \sigma(J_k(u)) \ = \ \infty.
\end{equation*}
Consequently, for every $\eta > 0$, there are infinitely many $k \in \N$ such that
\begin{equation*}
  \sigma(J_k(u)) \ \geq \ S_k^\eta.
\end{equation*}
Choose parameters $\beta$ and $\eta$ such that
\begin{equation*}
  0 \ < \ \eta \ < \ \delta\beta \quad \text{and} \quad 0 \ < \ \beta \ < \ \frac{(u-s)(1-\theta)}{\theta}.
\end{equation*}
Fix a sufficiently large $k \in \N$ such that $\sigma(J_k(u)) \geq S_k^\eta$, and write $R = S_k$.
Let $\mathcal{E} \subset J_k(u)$ be a maximal $R^\beta$-separated subset. Since the balls of radius $R^\beta$ centred at the points of $\mathcal{E}$ cover $J_k(u)$, the Frostman bound gives
\begin{equation*}
  \# \mathcal{E} \ \gtrsim \ \sigma(J_k(u)) R^{-\delta\beta} \ \geq \ R^{-(\delta\beta-\eta)}.
\end{equation*}
Since $G(2,1)$ is one-dimensional, we also have $\# \mathcal{E} \lesssim R^{-\beta}$. For each $V \in \mathcal{E}$, choose $x_V \in P_V X$ such that
\begin{equation*}
  N_{R^{1/\theta}} \big( P_V X \cap B(x_V,R) \big) \ > \ C_\theta^{-1} \biggl(\frac{R}{R^{1/\theta}}\biggr)^u.
\end{equation*}
Now define $r = (2R)^{1/\theta}$, the tubes $T_V$, and the sets $F_V$ and $E_V$ exactly as in the proof of Theorem \ref{thm:tube-lower-spectrum}, but with $u$ in place of $t$. The rest of the tube-overlap argument is unchanged: the lower bound for $\# \mathcal{E}$ above replaces the estimate $\# \mathcal{E} \gtrsim R^{-(\beta-\eta)}$, while the upper bound $\# \mathcal{E} \lesssim R^{-\beta}$ and the condition $\beta < (u-s)(1-\theta)/\theta$ play exactly the same roles as before. Consequently, one obtains a $2r$-separated subset of $X$ with cardinality at least
\begin{equation*}
  \gtrsim \ R^{-(\delta\beta-\eta)} \biggl(\frac{R}{R^{1/\theta}}\biggr)^u \ = \ (R^{-1/\theta})^{(1-\theta)u + (\delta\beta-\eta)\theta}.
\end{equation*}
Hence $\ubd X \geq (1-\theta)u + (\delta\beta-\eta)\theta$. Since this holds for every $\beta < (u-s)(1-\theta)/\theta$ and every $\eta \in (0,\delta\beta)$, letting $\beta \uparrow (u-s)(1-\theta)/\theta$ and $\eta \downarrow 0$ yields
\begin{equation*}
  \frac{\ubd X}{1-\theta} \ \geq \ u + \delta(u-s),
\end{equation*}
which contradicts \eqref{eq:delta-choice}. This proves \eqref{eq:tube-superlevel-bound}.

Now fix $\alpha \in [0,1)$. If
\begin{equation*}
  \as X \ = \ \frac{\ubd X}{1-\theta},
\end{equation*}
then
\begin{equation*}
  \frac{1}{1+\alpha} \biggl( \frac{\ubd X}{1-\theta} + \alpha \as X \biggr) \ = \ \frac{\ubd X}{1-\theta},
\end{equation*}
so the exceptional set is empty by the trivial bound $\as P_V X \leq \ubd X/(1-\theta)$ for all $V \in G(2,1)$. Thus assume that
\begin{equation*}
  \as X \ < \ \frac{\ubd X}{1-\theta}.
\end{equation*}
Set
\begin{equation*}
  u \ = \ \frac{1}{1+\alpha} \biggl( \frac{\ubd X}{1-\theta} + \alpha \as X \biggr).
\end{equation*}
Then $u \in (\as X,\ubd X/(1-\theta)]$, and a direct computation gives
\begin{equation*}
  \frac{\ubd X/(1-\theta) - u}{u - \as X} \ = \ \alpha.
\end{equation*}
Applying \eqref{eq:tube-superlevel-bound} with this choice of $u$ yields
\begin{equation*}
  \hd \{V \in G(2,1) : \as P_V X > u\} \ \leq \ \alpha,
\end{equation*}
which is the desired conclusion.

\section{Projection bounds for parametrised self-similar sets:~proof of Theorem \ref{thm:parametrised-spectrum}} \label{sec:parametrised-self-similar}

To prove Theorem \ref{thm:parametrised-spectrum}, we first show that the Assouad spectrum of $K_t$ coincides with the Assouad spectrum of the projection of $X$ onto a corresponding direction, and then apply the bounds from Theorems \ref{thm:tube-lower-spectrum} and \ref{thm:tube-lower-spectrum-alpha}. Fix $\theta \in (0,1)$ and $\alpha \in [0,1)$. Let $I$, the systems $\Psi_t$, the sets $K_t$, the planar lift $X$, and the number $s$ be as in the statement of the theorem. For each $t \in I$, define $\pi_t \colon \R^2 \to \R$ by
\begin{equation*}
  \pi_t(x,y) \ = \ x + t y.
\end{equation*}
Then, for each $i \in \{1,\ldots,N\}$ and $(x,y) \in \R^2$,
\begin{equation*}
  \pi_t(F_i(x,y)) \ = \ r x + a_i + t(r y + b_i) \ = \ r \pi_t(x,y) + a_i + t b_i \ = \ \psi_{i,t}(\pi_t(x,y)).
\end{equation*}
Therefore $\pi_t(X)$ is a non-empty compact set satisfying
\begin{equation*}
  \pi_t(X) \ = \ \bigcup_{i = 1}^N \psi_{i,t}(\pi_t(X)).
\end{equation*}
By uniqueness of the attractor of $\Psi_t$, this yields
\begin{equation} \label{eq:parametrised-kt}
  K_t \ = \ \pi_t(X).
\end{equation}
If $e_t = (1,t)/\sqrt{1+t^2} \in S^1$ and $V_t = \textup{span}(e_t) \in G(2,1)$, then $P_{V_t}(x,y) = (x,y) \cdot e_t = (x+t y)/\sqrt{1+t^2}$, and hence
\begin{equation} \label{eq:parametrised-scaling}
  \pi_t \ = \ \sqrt{1+t^2}\, P_{V_t}.
\end{equation}
Since the Assouad spectrum is bi-Lipschitz invariant, \eqref{eq:parametrised-kt} and \eqref{eq:parametrised-scaling} imply that
\begin{equation} \label{eq:parametrised-spectrum-equality}
  \as K_t \ = \ \as P_{V_t} X
\end{equation}
for all $t \in I$.

Let $\Gamma \colon I \to G(2,1)$ be given by $\Gamma(t) = V_t$. If $\sphericalangle(V_t,V_u)$ denotes the acute angle between $V_t$ and $V_u$, then
\begin{equation*}
  \sin \sphericalangle(V_t,V_u) \ = \ \sqrt{1-\frac{(1+tu)^2}{(1+t^2)(1+u^2)}} \ = \ \frac{|t-u|}{\sqrt{(1+t^2)(1+u^2)}}
\end{equation*}
for all $t,u \in I$. Since $\sphericalangle(V_t,V_u) \in [0,\pi/2]$, we have $\sin \sphericalangle(V_t,V_u) \approx \sphericalangle(V_t,V_u)$. As $I$ is bounded, the factor $((1+t^2)(1+u^2))^{-1/2}$ is bounded above and below by positive constants on $I \times I$, and so $\Gamma$ is bi-Lipschitz onto its image, and in particular injective. Consequently,
\begin{equation} \label{eq:parametrised-hd}
  \hd E \ = \ \hd \Gamma(E)
\end{equation}
for every $E \subseteq I$. Moreover, by the analyticity argument in the proof of Theorem \ref
{thm:tube-lower-spectrum-alpha}, the superlevel sets of the form $\{V \in G(2,1) : \as P_V X > c\}$ 
are measurable. Since $\gamma_{2,1}$ is a constant multiple of arc-length measure on $G(2,1)$, every measurable set $A \subseteq \Gamma(I)$ with $\gamma_{2,1}(A) = 0$ satisfies $\mathcal{H}^1(A) = 0$, where $\mathcal{H}^1$ is Hausdorff measure of dimension $1$. Because $\Gamma^{-1}$ is Lipschitz on $\Gamma(I)$, we also have $\mathcal{H}^1(\Gamma^{-1}(A)) = 0$. Since $\mathcal{H}^1$ agrees with Lebesgue measure $\mathcal{L}^1$ on $\R$, this yields
\begin{equation} \label{eq:parametrised-measure}
  \L^1(\Gamma^{-1}(A)) \ = \ 0
\end{equation}
whenever $A \subseteq \Gamma(I)$ is measurable and $\gamma_{2,1}(A) = 0$.

Let
\begin{equation*}
  u \ = \ \frac{1}{1+\alpha} \biggl( \frac{s}{1-\theta} + \alpha s \biggr).
\end{equation*}
Since $X$ is compact, the standard inequalities $\ubd X \leq \as X \leq \qad X$ recalled in Section \ref{intro} and the hypothesis $\qad X \leq s$ yield
\begin{equation*}
  \frac{1}{1+\alpha} \biggl( \frac{\ubd X}{1-\theta} + \alpha \as X \biggr) \ \leq \ u.
\end{equation*}
Therefore
\begin{equation*}
  \{V \in G(2,1) : \as P_V X > u\} \ \subseteq \ \biggl\{V \in G(2,1) : \as P_V X > \frac{1}{1+\alpha} \biggl( \frac{\ubd X}{1-\theta} + \alpha \as X \biggr)\biggr\}.
\end{equation*}
Hence Theorem \ref{thm:tube-lower-spectrum-alpha} gives
\begin{equation} \label{eq:parametrised-exceptional-directions}
  \hd \{V \in G(2,1) : \as P_V X > u\} \ \leq \ \alpha.
\end{equation}
By \eqref{eq:parametrised-spectrum-equality}, we have $\Gamma (\{t \in I : \as K_t > u\}) \subseteq \{V \in G(2,1) : \as P_V X > u\}$. Combining this inclusion with \eqref{eq:parametrised-hd} and \eqref{eq:parametrised-exceptional-directions}, we obtain
\begin{equation*}
  \hd \{t \in I : \as K_t > u\} \ = \ \hd \Gamma (\{t \in I : \as K_t > u\}) \ \leq \ \alpha.
\end{equation*}
This proves the exceptional-set estimate.

For the almost-sure estimate, Theorem \ref{thm:tube-lower-spectrum} and \eqref{eq:parametrised-dim-chain} imply that
\begin{equation*}
  \as P_V X \ \leq \ \frac{1}{2} \biggl( \frac{\ubd X}{1-\theta} + \as X \biggr) \ \leq \ \frac{1}{2} \biggl( \frac{s}{1-\theta} + s \biggr)
\end{equation*}
for $\gamma_{2,1}$-almost all $V \in G(2,1)$. Hence the set
\begin{equation*}
  A \ = \ \biggl\{V \in \Gamma(I) : \as P_V X > \frac{1}{2} \biggl( \frac{s}{1-\theta} + s \biggr)\biggr\}
\end{equation*}
has $\gamma_{2,1}$-measure zero. By \eqref{eq:parametrised-spectrum-equality},
\begin{equation*}
  \biggl\{t \in I : \as K_t > \frac{1}{2} \biggl( \frac{s}{1-\theta} + s \biggr)\biggr\} \ = \ \Gamma^{-1}(A).
\end{equation*}
Hence \eqref{eq:parametrised-measure} gives
\begin{equation*}
  \L^1 \biggl(\biggl\{t \in I : \as K_t > \frac{1}{2} \biggl( \frac{s}{1-\theta} + s \biggr)\biggr\}\biggr) \ = \ 0
\end{equation*}
as claimed.

\section{Projections with two distinct Assouad spectrum values:~proof of Theorem \ref{main1}} \label{sec:projection-two-values}

To prove Theorem \ref{main1}, we build a planar set from thin, widely separated pieces (which we call `bricks')  whose orientations alternate between two directions. Each piece is shaped so that projecting it along its thin side reveals a set of prescribed dimension, while projecting along the thick side collapses it to something negligible. Directions close to the first coordinate axis see one family of pieces and pick up dimension $s$, while directions close to the second coordinate axis see the other family and pick up dimension $t$.

\subsection{Construction of the set} \label{sec:construction-of-the-set}

The set $X$ is constructed as a disjoint union of sparse bricks, alternating in orientation, so that projections onto different directions capture the dimensions of two distinct families of bricks. We begin by choosing sequences $(\alpha_n)_{n=1}^\infty$ in $(0,1)$ and $(\beta_n)_{n=1}^\infty$ in $[1,\infty)$ such that for all $\theta \in (0,1)$,
\begin{equation} \label{alpha0}
  \lim_{n \to \infty} \alpha_n \ = \ 0, \qquad \lim_{n \to \infty} \alpha_n^{1/\theta} n \ = \ \infty
\end{equation}
and
\begin{equation} \label{betasmall}
  \lim_{n \to \infty} \beta_n \ = \ \infty, \qquad \lim_{n \to \infty} \frac{\beta_n}{\alpha_n^{1/\theta} n} \ = \ 0.
\end{equation}
For concreteness, we take $\alpha_n = (\log(n+2))^{-1}$ and $\beta_n = n^{1/2}$ for all $n \geq 1$. To verify \eqref{alpha0} and \eqref{betasmall} for all $\theta \in (0,1)$, note that $\alpha_n^{1/\theta} = (\log(n+2))^{-1/\theta}$, so $\alpha_n^{1/\theta} n = n/(\log(n+2))^{1/\theta}$. Since $(\log(n+2))^{1/\theta}$ grows slower than any positive power of $n$, we have $\alpha_n^{1/\theta} n \to \infty$ for all $\theta \in (0,1)$. Moreover,
\begin{equation*}
  \frac{\beta_n}{\alpha_n^{1/\theta} n} \ = \ \frac{(\log(n+2))^{1/\theta}}{n^{1/2}} \to 0
\end{equation*}
for all $\theta \in (0,1)$. We also set
\begin{equation*}
  \delta_n \ = \ (\log(n+2))^{-1/2}.
\end{equation*}
Then $\delta_n \to 0$ as $n \to \infty$. Fix $\theta \in (0,1)$ and $s, t, u \in (0,1]$. Let
\begin{equation*}
  E_{u}^n \ = \ \biggl\{ \sum_{k = 1}^N \pm 2^{-k/{u}}\biggr\} \quad\text{and}\quad F_{u}^n \ = \ \biggl\{ \sum_{k = 1}^N \pm 2^{-k/\delta_n}\biggr\}
\end{equation*}
where $N = N(n,u)$ is chosen such that
\begin{equation} \label{defN}
  2^{N-1} \ \leq \ \alpha_n^{u(1-1/\theta)} \ \leq \ 2^{N}.
\end{equation}
The set $E_u^n \subseteq [-1,1]$ is a finite set of $2^N$ points approximating a self-similar set of dimension $u$ at scale $2^{-N/{u}} \approx \alpha_n^{1/\theta-1}$. Similarly, $F_u^n \subseteq [-1,1]$ is a finite set of $2^N$ points approximating a self-similar set of dimension $\delta_n$ at scale $2^{-N/\delta_n}$.

We establish the covering estimates needed below. Let $q_u = 2^{-1/u} \leq 1/2$, and let $K_u$ be the self-similar set generated by the maps $x \mapsto q_u x \pm q_u$. If $u < 1$, then $q_u < 1/2$, these maps satisfy the strong separation condition, and the natural Bernoulli measure on $K_u$ is $u$-Ahlfors regular. Consequently, if an interval $J$ has length at least $q_u^{m+1}$, then the number of level-$m$ cylinders meeting $J$ is at most $\lesssim |J|^u 2^m$, uniformly in $m \in \mathbb{N}$: indeed, the union of these cylinders is contained in a $Cq_u^m$-neighborhood of $J$, which still has length $\lesssim |J|$, and each level-$m$ cylinder has Bernoulli measure $2^{-m}$. Let $E$ be a scaled copy of $E_u^n$ of diameter $L$, let $J$ be an interval of length $\ell \leq L$, and choose $m$ such that $Lq_u^{m+1} < \rho \leq Lq_u^m$. Since $\rho \leq \ell$, the normalized interval has length $\ell/L \geq q_u^{m+1}$, and therefore $E \cap J$ is contained in the union of at most $\lesssim (\ell/L)^u 2^m \lesssim (\ell/\rho)^u$ level-$m$ cylinders. After normalizing by $L$, the points of $E \cap J$ lying in a fixed level-$m$ cylinder differ only in the tail $\sum_{k = m+1}^N \pm q_u^k$, so each such cluster has diameter $\lesssim Lq_u^{m+1} \lesssim \rho$. After subdividing by an absolute constant if necessary, each such cluster can be covered by at most $\lesssim 1$ intervals of length $\rho$. Thus every scaled copy of $E_u^n$ of diameter $L$ can be covered inside any interval of length $\ell \leq L$ by at most $\lesssim (\ell/\rho)^u$ intervals of length $\rho$, uniformly in $0 < \rho \leq \ell$ when $u < 1$. If $u = 1$, then $E_1^n$ is an arithmetic progression with gap comparable to $2^{-N}$. Hence a scaled copy of $E_1^n$ of diameter $L$ can be covered inside any interval of length $\ell \leq L$ by at most $\lesssim \ell/\rho$ intervals of length $\rho$, uniformly in $0 < \rho \leq \ell$, so the same covering estimate also holds at the endpoint $u = 1$. The same argument with $q_n = 2^{-1/\delta_n}$ shows that every scaled copy of $F_u^n$ of diameter $L$ can be covered inside any interval of length $\ell \leq L$ by at most $\lesssim (\ell/\rho)^{\delta_n}$ intervals of length $\rho$, uniformly in $0 < \rho \leq \ell$.

Since $N \lesssim \log\log n$ uniformly in $u \in (0,1]$ and $\delta_n^{-1} = (\log(n+2))^{1/2}$, we have $N/\delta_n = o(\log n)$. Hence $2^{-N/\delta_n} = n^{-o(1)}$, and therefore
\begin{equation} \label{longgap}
  \frac{\alpha_n}{\beta_n 2^{-N/\delta_n}} \ \to \ 0
\end{equation}
uniformly in $u \in (0,1]$. Let
\begin{equation*}
  Q(n,u) \ = \ \{ (y_k, z_k) \}_{k = 1}^{2^N} \ \subset \ [0,\alpha_n] \times [0,\beta_n]
\end{equation*}
be a finite set consisting of $2^N$ points lying on the graph of a strictly increasing function. We require that the set of first coordinates, $\{y_k\}_{k=1}^{2^N}$, is a similar copy of $E_u^n$ scaled to have diameter $\alpha_n$ and the set of second coordinates, $\{z_k\}_{k=1}^{2^N}$, is a similar copy of $F_u^n$ scaled to have diameter $\beta_n$. To realise such a set, choose these two scaled copies, order them increasingly as $0 \leq y_1 < \cdots < y_{2^N} \leq \alpha_n$ and $0 \leq z_1 < \cdots < z_{2^N} \leq \beta_n$, and pair points with the same rank. Then $Q(n,u)$ lies on the graph of a strictly increasing function. Translations and rotations of the sets $Q(n,u)$ are called \emph{bricks}. Let $(x_n)_{n=1}^\infty$ be a sequence of translations in $\mathbb{R}^2$ defined recursively by
\begin{equation*}
  x_1 \ = \ (0,0), \qquad x_{n+1} \ = \ x_n + (4\beta_n+4, 4\beta_n+4).
\end{equation*}
We equip $G(2,1)$ with the angular metric: for $V = \textup{span}(e^{i\phi})$ and $W = \textup{span}(e^{i\psi})$ with $\phi, \psi \in [0,\pi)$, we set $\|V - W\|$ equal to the acute angle between $V$ and $W$, that is, $\|V - W\| = \min\{|\phi - \psi|, \pi - |\phi - \psi|\}$. Let $I_s \subseteq G(2,1)$ be the ball of radius $1/100$ centred at the first coordinate axis and let $I_t \subseteq G(2,1)$ be the ball of radius $1/100$ centred at the second coordinate axis. The radius $1/100$ is chosen sufficiently small to ensure disjointness. Let $\{V_n^s\}$ be a countable dense subset of $I_s$ such that for all $V \in I_s$ there are infinitely many $n \in \mathbb{N}$ satisfying
\begin{equation*}
  \| V - V_n^s \| \ \lesssim \ 1/n.
\end{equation*}
Such a set $\{V_n^s\}$ can be constructed as follows. For each $k \in \mathbb{N}$, let $A_k^s \subseteq I_s$ be a maximal $2^{-k}$-separated set. Since $I_s$ is an interval in the angular metric, the set $A_k^s$ is also a $2^{-k}$-net of $I_s$ and satisfies $\# A_k^s \approx 2^k$. Enumerate the points of the sets $A_k^s$ stage by stage to obtain a sequence $\{V_n^s\}_{n=1}^\infty$. If $V \in I_s$, then for every $k \in \mathbb{N}$ there exists a point of $A_k^s$ within distance $\lesssim 2^{-k}$ of $V$, and the corresponding index satisfies $n \approx 2^k$. Hence $\| V - V_n^s \| \lesssim 2^{-k} \lesssim 1/n$ for infinitely many $n \in \mathbb{N}$. Similarly, let $\{V_n^t\}$ be a countable dense subset of $I_t$ such that for all $V \in I_t$ there are infinitely many $n \in \mathbb{N}$ satisfying
\begin{equation*}
  \| V - V_n^t \| \ \lesssim \ 1/n.
\end{equation*}
For each $k \in \mathbb{N}$, let $A_k^t \subseteq I_t$ be a maximal $2^{-k}$-separated set, and enumerate the points of the sets $A_k^t$ stage by stage. The same argument yields the required approximation property for $\{V_n^t\}$. Finally, for each $n$, let $\phi_{2n}, \phi_{2n-1} \in [0,\pi)$ be such that $V_n^s = \textup{span}(e^{i \phi_{2n}})$ and $V_n^t = \textup{span}(e^{i \phi_{2n-1}})$.

We define
\begin{equation*}
  X \ = \ \bigcup_{n = 1}^\infty \Big(Q(2n,s) e^{i \phi_{2n}} +x_{2n} \Big) \ \cup \ \bigcup_{n = 1}^\infty \Big(Q(2n-1,t) e^{i \phi_{2n-1}} +x_{2n-1} \Big).
\end{equation*}
We claim that $X$ is uniformly discrete. Write $q_n = 2^{-1/\delta_n} < 1/2$. If two points of $F_u^n$ first differ in the $m$th sign, then their difference is $2q_n^m$ plus a tail of size $O(q_n^{m+1})$, so, since $q_n < 1/2$, their distance is comparable to $q_n^m$, and in particular the minimal gap between distinct points of $F_u^n$ is comparable to $q_n^N$, uniformly in $u \in (0,1]$. Since the long coordinates in $Q(n,u)$ form a scaled copy of $F_u^n$ of diameter $\beta_n$, it follows that the minimal gap between two distinct long coordinates in $Q(n,u)$ is comparable to $\beta_n 2^{-N/\delta_n}$. Euclidean distance dominates the difference of the long coordinates, and rotations preserve Euclidean distance, so the minimal distance between two distinct points of the rotated brick $Q(n,u)e^{i\phi}$ is also $\gtrsim \beta_n 2^{-N/\delta_n}$. Since $N/\delta_n = o(\log n)$ and $\beta_n = n^{1/2}$, we have $\beta_n 2^{-N/\delta_n} \to \infty$. Therefore, there exists $n_0 \in \mathbb{N}$ such that every brick with index at least $n_0$ is $1$-separated.

Moreover, every point of the rotated brick $Q(n,u)e^{i\phi}+x_n$ lies in $x_n + [-(\beta_n+1),\beta_n+1]^2$. If $m > n$, then the difference between the centres $x_m$ and $x_n$ in each coordinate is $\sum_{k = n}^{m-1} (4\beta_k+4)$. Since $(\beta_n)_{n=1}^\infty$ is increasing and $\beta_{k+1} \leq 2\beta_k$ for all $k$, this difference is at least $4\beta_{m-1}+4$, whereas
\begin{equation*}
  (\beta_n+1) + (\beta_m+1) \ \leq \ 3\beta_{m-1}+2.
\end{equation*}
Hence distinct bricks are separated by distance at least $1$. The finitely many points belonging to the first $n_0-1$ bricks also have a positive minimal separation, and therefore $X$ is uniformly discrete. Since each brick is non-empty and $|x_n| \to \infty$ by construction, the set $X$ is also unbounded. In particular, $X$ is countable and closed.

\subsection{Lower bound for the Assouad spectrum of projections} \label{sec:projection-two-values-lower}

Let $V \in I_s$. We show that $\as P_V X \geq s$. By construction we can find infinitely many $n \in \mathbb{N}$ such that
\begin{equation} \label{quantsep}
  \| V - V_n^s \| \ \lesssim \ 1/n.
\end{equation}
Fix such an $n$. The projection onto $V_n^s$ of the single brick
\begin{equation} \label{setdis}
  Q(2n,s) e^{i \phi_{2n}} +x_{2n}.
\end{equation}
consists of $2^N$ points which are $\alpha_{2n} 2^{-N/{s}} \approx \alpha_{2n}^{1/\theta}$ separated and lie in a common interval of length $\alpha_{2n}$ by \eqref{defN}. Let $(y_i,z_i)$ and $(y_j,z_j)$ be two points of $Q(2n,s)$ with $j > i$. Since $Q(2n,s)$ lies on the graph of a strictly increasing function, we have $y_j-y_i > 0$ and $z_j-z_i > 0$. Writing $\delta = \| V - V_n^s \|$, we obtain
\begin{equation*}
  |P_V((y_j,z_j)e^{i \phi_{2n}} + x_{2n}) - P_V((y_i,z_i)e^{i \phi_{2n}} + x_{2n})| \ \geq \ (y_j-y_i)\cos \delta - (z_j-z_i)|\sin \delta|.
\end{equation*}
Since $y_j-y_i \gtrsim \alpha_{2n} 2^{-N/{s}} \approx \alpha_{2n}^{1/\theta}$, $z_j-z_i \leq \beta_{2n}$, and $|\sin \delta| \lesssim \delta$, \eqref{quantsep} shows that the projection of \eqref{setdis} onto $V$ consists of $2^{N} \approx \alpha_{2n}^{s(1-1/\theta)}$ points which, by \eqref{betasmall}, are
\begin{equation*}
  \gtrsim \ \alpha_{2n}^{1/\theta} - \beta_{2n} \| V - V_n^s \| \ \gtrsim \ \alpha_{2n}^{1/\theta}
\end{equation*}
separated and lie in a common interval of length
\begin{equation*}
  \lesssim \ \alpha_{2n} + \beta_{2n} \| V - V_n^s \| \ \lesssim \ \alpha_{2n}.
\end{equation*}
Choose $y_n \in P_V(Q(2n,s)e^{i \phi_{2n}} + x_{2n}) \subseteq P_V X$. Since the projected points lie in an interval of length $\lesssim \alpha_{2n}$, there exists an absolute constant $C \geq 1$ such that
\begin{equation*}
  P_V(Q(2n,s)e^{i \phi_{2n}} + x_{2n}) \ \subseteq \ B(y_n,C\alpha_{2n})
\end{equation*}
for all sufficiently large $n$. Since these projected points lie on a line and are $\gtrsim \alpha_{2n}^{1/\theta}$ separated, every ball of radius $(C\alpha_{2n})^{1/\theta}$ meets at most $\lesssim 1$ of them, where the implied constant depends only on $C$ and $\theta$. By monotonicity of the covering number in the set being covered, this yields
\begin{equation*}
  N_{(C\alpha_{2n})^{1/\theta}}(P_V X \cap B(y_n,C\alpha_{2n})) \ \geq \ N_{(C\alpha_{2n})^{1/\theta}}(P_V(Q(2n,s)e^{i \phi_{2n}} + x_{2n})) \ \gtrsim \ 2^N \ \approx \ \alpha_{2n}^{s(1-1/\theta)}.
\end{equation*}
If $\sigma < s$, then
\begin{equation*}
  \frac{N_{(C\alpha_{2n})^{1/\theta}}(P_V X \cap B(y_n,C\alpha_{2n}))}{(C\alpha_{2n})^{\sigma(1-1/\theta)}} \ \approx \ \alpha_{2n}^{(s-\sigma)(1-1/\theta)} \to \infty
\end{equation*}
since $1 - 1/\theta < 0$ and $\alpha_{2n} \to 0$ by \eqref{alpha0}. Thus no exponent $\sigma < s$ can satisfy the definition of $\as P_V X$. Since the estimate holds for infinitely many $n$, we conclude that $\as P_V X \geq s$, as required.

The same argument gives $\as P_V X \geq t$ for all $V \in I_t$.

\subsection{Upper bound for the Assouad spectrum of projections} \label{sec:projection-two-values-upper}

Let $V \in I_s$. We show that $\as P_V X \leq s$. Fix any $x \in P_V X$, set $r = R^{1/\theta}$ for $R \in (0,1)$, and consider $P_V X \cap B(x,R)$. Since $\beta_n = n^{1/2}$, we have $\beta_{n+1} \leq 2\beta_n$ for all $n$. Let $a_0 = \cos(1/100)-\sin(1/100) > 3/4$ and choose $c_0 \in (1,4a_0/3)$. Since the balls $I_s$ and $I_t$ have radius $1/100$ and $\alpha_n/\beta_n \to 0$, there exists $n_1 \in \mathbb{N}$ such that, for every $n \geq n_1$, the projection of the $n$th brick onto $V$ is contained in an interval centred at $P_V x_n$ of radius at most $c_0\beta_n$: if the short axis of the brick lies in $I_s$, then the projection radius is at most $\alpha_n + \beta_n\sin(1/50)$, whereas if the short axis lies in $I_t$, then the projection radius is at most $\alpha_n + \beta_n$. Moreover, every direction in $I_s$ makes angle at most $1/100$ with the first coordinate axis, so the translation vector $(4\beta_n+4,4\beta_n+4)$ has projection at least $a_0(4\beta_n+4)$ onto $V$. Consequently, consecutive projected bricks with indices $n$ and $n+1$ are separated by distance at least
\begin{equation*}
  a_0(4\beta_n+4) - c_0\beta_n - c_0\beta_{n+1} \ \geq \ (4a_0-3c_0)\beta_n + 4a_0 \ > \ 2
\end{equation*}
for all $n \geq n_1$. Therefore every ball $B(x,R)$ with $R \in (0,1)$ meets the projection of at most one brick with index at least $n_1$, and it suffices to treat one such brick, since the finitely many bricks with index below $n_1$ contribute at most a bounded constant factor.

First, suppose that the intersecting brick is an even brick, $Q(2n,s)e^{i\phi_{2n}} + x_{2n}$, and identify $Q(2n,s)$ with a subset of the $uv$-plane where the $u$-axis is the short side and the $v$-axis is the long side. Since $V$ and $V_n^s$ both lie in $I_s$, there exist absolute constants $0 < c_2 \leq C_2 < \infty$ and numbers $a_n,b_n,\lambda_n \in \R$ with $c_2 \leq |a_n| \leq C_2$ and $|\lambda_n| \leq C_2$ such that
\begin{equation*}
  P_V((u,v)e^{i\phi_{2n}} + x_{2n}) \ = \ P_V x_{2n} + b_n + a_n(u + \lambda_n v)
\end{equation*}
for all $(u,v) \in Q(2n,s)$. Indeed, in these coordinates the projection is an affine linear functional, and because both $V$ and $V_n^s$ lie in the fixed interval $I_s$, the coefficient of the short coordinate stays bounded away from zero while all coefficients remain uniformly bounded. Passing from the actual projection to the model map $(u,v) \mapsto u + \lambda_n v$ therefore only rescales lengths by an absolute factor and adds a translation. More precisely, if $A \subseteq \R$ is bounded and $\rho > 0$, then
\begin{equation*}
  N_{\rho/C_2}(A) \ \leq \ N_\rho(a_n A + b_n) \ \leq \ N_{\rho/c_2}(A).
\end{equation*}
Thus an estimate for the model map $(u,v) \mapsto u + \lambda_n v$ at scales comparable to $R$ and $r$ transfers back to the actual projection with only an absolute multiplicative constant. We claim that
\begin{equation*}
  N_r(P_V X \cap B(x,R)) \ \lesssim \ \biggl(\frac{R}{r}\biggr)^{s+\delta_{2n}}.
\end{equation*}
If $\alpha_{2n} \geq R$, then the relevant points can be covered individually. By \eqref{defN},
\begin{equation*}
  N_r(P_V X \cap B(x,R)) \ \leq \ 2^N \ \approx \ \alpha_{2n}^{s(1-1/\theta)} \ \lesssim \ R^{s(1-1/\theta)} \ = \ \biggl(\frac{R}{r}\biggr)^s.
\end{equation*}
Hence it remains to consider the case $\alpha_{2n} < R$. If $\lambda_n = 0$, then the projection of the brick is, up to translation and bounded scaling, a scaled copy of the first coordinate set. If $\alpha_{2n} < r$, then the projected points can be covered by $\lesssim 1$ balls of radius $r$. Otherwise, since the first coordinates form a scaled copy of $E_s^{2n}$ of diameter $\alpha_{2n}$, the covering estimate from Section~\ref{sec:construction-of-the-set} gives at most
\begin{equation*}
  \lesssim \ \biggl(\frac{\alpha_{2n}}{r}\biggr)^s \ \leq \ \biggl(\frac{R}{r}\biggr)^s
\end{equation*}
balls. Assume now that $\lambda_n \neq 0$. The relevant points of $Q(2n,s)$ satisfy
\begin{equation*}
  |u + \lambda_n v - c| \ \leq \ CR
\end{equation*}
for some $c \in \R$. Since $0 \leq u \leq \alpha_{2n} < R$, it follows that $\lambda_n v$ lies in an interval of length $\lesssim R$. Hence the corresponding long coordinates lie in an interval $J$ of length $\lesssim R/|\lambda_n|$. We cover the intersection of the long coordinates with $J$ by intervals of length $\min\{r/|\lambda_n|,\beta_{2n}\}$. Since the long coordinates form a scaled copy of $F_s^{2n}$ with dimension $\delta_{2n}$, the covering estimate from Section~\ref{sec:construction-of-the-set} gives
\begin{equation*}
  N_v \ \lesssim \ \biggl(\frac{\min\{R/|\lambda_n|,\beta_{2n}\}}{\min\{r/|\lambda_n|,\beta_{2n}\}}\biggr)^{\delta_{2n}}.
\end{equation*}
If $|\lambda_n|\beta_{2n} \leq r$, then both minima equal $\beta_{2n}$ and $N_v \lesssim 1$. If $r < |\lambda_n|\beta_{2n} \leq R$, then the numerator equals $\beta_{2n}$ and the denominator equals $r/|\lambda_n|$, so
\begin{equation*}
  N_v \ \lesssim \ \biggl(\frac{|\lambda_n|\beta_{2n}}{r}\biggr)^{\delta_{2n}} \ \leq \ \biggl(\frac{R}{r}\biggr)^{\delta_{2n}}.
\end{equation*}
If $|\lambda_n|\beta_{2n} > R$, then the numerator equals $R/|\lambda_n|$ and the denominator equals $r/|\lambda_n|$, so
\begin{equation*}
  N_v \ \lesssim \ \biggl(\frac{R}{r}\biggr)^{\delta_{2n}}.
\end{equation*}
Thus in every case
\begin{equation*}
  N_v \ \lesssim \ \biggl(\frac{R}{r}\biggr)^{\delta_{2n}}.
\end{equation*}
For each interval $W$ in this cover, the variation of $\lambda_n v$ over $W$ is at most $r$. Therefore the projection of the corresponding points is contained in the $r$-neighborhood of a translate of the associated first coordinates. If $\alpha_{2n} < r$, then these projected points can be covered by $\lesssim 1$ balls of radius $r$. Otherwise, since the first coordinates form a scaled copy of $E_s^{2n}$ of diameter $\alpha_{2n}$, the covering estimate from Section~\ref{sec:construction-of-the-set} gives at most
\begin{equation*}
  \lesssim \ \biggl(\frac{\alpha_{2n}}{r}\biggr)^s \ \leq \ \biggl(\frac{R}{r}\biggr)^s
\end{equation*}
balls. Hence each interval $W$ contributes at most $\lesssim (R/r)^s$ balls, and therefore
\begin{equation*}
  N_{r}(P_V X \cap B(x,R)) \ \lesssim \ \biggl(\frac{R}{r}\biggr)^{s+\delta_{2n}}.
\end{equation*}

Next, suppose that the intersecting brick is an odd brick, $Q(2n-1,t)e^{i\phi_{2n-1}} + x_{2n-1}$, and identify $Q(2n-1,t)$ with a subset of the $uv$-plane where the $u$-axis is the short side and the $v$-axis is the long side. Since $V \in I_s$ while $V_n^t \in I_t$, after translation and bounded scaling the projection onto $V$ is of the form
\begin{equation*}
  (u,v) \ \mapsto \ \eta_n u + \lambda_n v
\end{equation*}
with $|\eta_n| \lesssim 1$ and $|\lambda_n| \approx 1$, where the implied constants are uniform in $n$ because directions in $I_s$ and $I_t$ stay a fixed positive angle apart. Let $(u_i,v_i)$ and $(u_j,v_j)$ be two points of $Q(2n-1,t)$ with $v_i < v_j$. Since $Q(2n-1,t)$ lies on the graph of a strictly increasing function, we also have $u_i < u_j$, and therefore
\begin{equation*}
  |\lambda_n||v_j-v_i| - |\eta_n||u_j-u_i| \ \leq \ |\eta_n(u_j-u_i) + \lambda_n(v_j-v_i)| \ \leq \ |\lambda_n||v_j-v_i| + |\eta_n||u_j-u_i|.
\end{equation*}
Since $0 \leq u \leq \alpha_{2n-1}$ throughout the brick, we have $|u_j-u_i| \leq \alpha_{2n-1}$. Moreover, because the long coordinates form a scaled copy of $F_t^{2n-1}$ of diameter $\beta_{2n-1}$, every distinct pair satisfies $|v_j-v_i| \gtrsim \beta_{2n-1} 2^{-N/\delta_{2n-1}}$. By \eqref{longgap}, for all sufficiently large $n$ and all distinct pairs of points in $Q(2n-1,t)$ we therefore have
\begin{equation*}
  |\eta_n||u_j-u_i| \ \leq \ \frac12 |\lambda_n||v_j-v_i|.
\end{equation*}
Hence
\begin{equation*}
  \frac12 |\lambda_n||v_j-v_i| \ \leq \ |\eta_n(u_j-u_i) + \lambda_n(v_j-v_i)| \ \leq \ 2 |\lambda_n||v_j-v_i|.
\end{equation*}
Since $|\lambda_n| \approx 1$, it follows that, for large $n$, the projected odd brick is bi-Lipschitz to a scaled copy of the long coordinate set $F_t^{2n-1}$, with bi-Lipschitz constants independent of $n$. Therefore every interval of length $R$ in the projection corresponds to an interval of long coordinates of length $\lesssim R$, and intervals of long-coordinate radius comparable to $r$ project to intervals of radius comparable to $r$. Replacing both $R$ and $r$ by comparable multiples changes covering numbers by at most an absolute constant. Since scaled copies of $F_t^{2n-1}$ satisfy the covering estimate stated in Section~\ref{sec:construction-of-the-set}, every such projected interval can be covered at scale $r$ by at most $\lesssim (\frac{R}{r})^{\delta_{2n-1}}$ balls. Consequently,
\begin{equation*}
  N_r(P_V X \cap B(x,R)) \ \lesssim \ \biggl(\frac{R}{r}\biggr)^{\delta_{2n-1}}
\end{equation*}
for all sufficiently large odd indices.

Fix $\varepsilon > 0$. Since $\delta_n \to 0$, we may choose $n_0 \geq n_1$ so large that $\delta_n < \varepsilon$ for all $n \geq n_0$. The preceding estimates then show that every brick with index at least $n_0$ contributes at most
\begin{equation*}
  \lesssim \ \biggl(\frac{R}{r}\biggr)^{s+\varepsilon}
\end{equation*}
covering balls. On the other hand, there are only finitely many bricks with index smaller than $n_0$, and each of them contains only finitely many points, so their contributions can be absorbed into the multiplicative constant. We conclude that
\begin{equation*}
  N_{R^{1/\theta}}(P_V X \cap B(x,R)) \ \lesssim_\varepsilon \ \biggl(\frac{R}{R^{1/\theta}}\biggr)^{s+\varepsilon}
\end{equation*}
for all $x \in P_V X$ and all $R \in (0,1)$. Since $\varepsilon > 0$ was arbitrary, this yields $\as P_V X \leq s$.

By the same argument, with the roles of the even and odd bricks interchanged, we obtain $\as P_V X \leq t$ for all $V \in I_t$.

\section{Infinitely many distinct values:~proof of Theorem \ref{main1+}} \label{sec:projection-infinite-values}

To prove Theorem \ref{main1+}, we keep the same global sequences $(\alpha_n)_{n=1}^{\infty}$, $(\beta_n)_{n=1}^{\infty}$, and $\delta_n$ as in Section~\ref{sec:construction-of-the-set}. The geometric inputs that we reuse from the proof of Theorem~\ref{main1} are of two kinds, according to the orientation of the short axis of a brick relative to a direction $V$. Suppose first that the short axis approximates $V$ at rate $\lesssim 1/n$. Then the projection onto $V$ contains a witness ball of radius $C\alpha_n$ carrying $2^N$ points that are $\alpha_n^{1/\theta}$-separated, where $N$ is given by \eqref{defN}; moreover, if the brick is built from the exponent $u$, the same-family estimate of Section~\ref{sec:projection-two-values-upper} covers every interval of length $R$ in such a projection, at scale $R^{1/\theta}$, by $\lesssim (R/R^{1/\theta})^{u+\delta_n}$ balls. Suppose instead that the short axis makes a fixed positive angle with $V$. Then \eqref{longgap} forces the projection to be bi-Lipschitz to the long coordinate set, so it contributes only the exponent $\delta_n$.

For Theorem \ref{main1+}, we are given a countable family of exponents $s_1,s_2,\ldots\in (0,1]$
and we want a single closed set $X\subset \mathbb{R}^2$ together with pairwise disjoint open sets $I_j\subset G(2,1)$, $j\in \N$, such that for all $j\in \N$ and all $V\in I_j$, we have
\begin{equation*}
  \as P_V X \ = \ s_j.
\end{equation*}
The proof uses countably many direction families, one for each prescribed value, but only one family at each dyadic scale.

\subsection{Choice of directions in the Grassmannian}

For each $j \in \N$, let $W_j = \textup{span}(e^{i2^{-j-10}})$ and define
\begin{equation*}
  I_j \ = \ B(W_j,2^{-j-13}).
\end{equation*}
Then each $I_j$ is a non-empty open ball in $G(2,1)$, all the balls are contained in the ball of radius $1/100$ centred at the first coordinate axis, and the closures are pairwise disjoint. Since $I_j$ is an interval of length $2^{-j-12}$ in the angular metric, every $2^{-k}$-separated subset of $I_j$ has cardinality at most $1 + 2^{-j-12}2^k \leq 2^k$, and therefore every maximal $2^{-k}$-separated subset of $I_j$ has cardinality at most $2^k$. For each fixed $j$, let
\begin{equation*}
  d_j \ = \ \dist\biggl(I_j,\bigcup_{l \neq j} I_l\biggr) \ > \ 0.
\end{equation*}
Then every $V \in I_j$ has distance at least $d_j$ from $\bigcup_{l \neq j} I_l$.

\subsection{Indexing the bricks}

We keep the bricks $Q(n,u)$ from Section~\ref{sec:construction-of-the-set}: for each $u\in(0,1]$ and $n\in \N$, the set $Q(n,u)$ lies in $[0,\alpha_n]\times[0,\beta_n]$ and is arranged on the graph of an increasing function so that its coordinate projections are scaled copies of $E_u^n$ and $F_u^n$. Fix a partition of $\mathbb{N}$ into infinite subsets $K_1,K_2,\ldots$. For each $j$ and each $k \in K_j$, let $\{V_{j,k,\ell}\}_{\ell=1}^{M_{j,k}} \subset I_j$ be a maximal $2^{-k}$-separated subset of $I_j$. By the choice of $I_j$, we have $M_{j,k} \le 2^k$. Choose distinct indices
\begin{equation*}
  n_{j,k,\ell} \ \in \ \{2^k,\ldots,2^{k+1}-1\}
\end{equation*}
for $\ell \in \{1,\ldots,M_{j,k}\}$, which is possible because the dyadic block contains $2^k$ indices. For each $\ell$, choose $\phi_{j,k,\ell} \in [0,\pi)$ such that
\begin{equation*}
  V_{j,k,\ell} \ = \ \textup{span}(e^{i\phi_{j,k,\ell}}),
\end{equation*}
and define
\begin{equation*}
  Q_{j,k,\ell} \ = \ Q(n_{j,k,\ell},s_j)e^{i\phi_{j,k,\ell}} + x_{n_{j,k,\ell}}.
\end{equation*}
For every fixed $V \in I_j$ and every $k \in K_j$, maximality gives an index $\ell(k)$ such that
\begin{equation*}
  \|V - V_{j,k,\ell(k)}\| \ \lesssim \ 2^{-k} \ \lesssim \ \frac{1}{n_{j,k,\ell(k)}}.
\end{equation*}
Finally, we define
\begin{equation*}
  X \ = \ \bigcup_{j = 1}^\infty \ \bigcup_{k \in K_j} \ \bigcup_{\ell = 1}^{M_{j,k}} Q_{j,k,\ell}.
\end{equation*}
Since the sets $K_j$ are disjoint and the dyadic blocks $\{2^k,\ldots,2^{k+1}-1\}$ are pairwise disjoint, each index $n$ is used for at most one brick. The separation argument from Section~\ref{sec:construction-of-the-set} depends only on this fact and on the translations $x_n$, so it applies verbatim and shows that $X$ is uniformly discrete. Hence $X$ is locally finite, and therefore closed. Moreover, if $k \in K_j$, then the maximal $2^{-k}$-separated subset of $I_j$ is non-empty, and hence there is at least one brick with index in $\{2^k,\ldots,2^{k+1}-1\}$. Since $|x_n| \to \infty$ with $n$, it follows that $X$ is unbounded.

\subsection{Bounding the Assouad spectrum of projections}

Fix $j\in\mathbb{N}$ and $V\in I_j$. For every $k \in K_j$, choose $\ell(k)$ as above and set $n_k = n_{j,k,\ell(k)}$. Since $n_k \in [2^k,2^{k+1})$, we have $n_k \to \infty$ as $k \to \infty$. By the argument from Section~\ref{sec:projection-two-values-lower}, the set $P_V Q_{j,k,\ell(k)}$ is contained in a ball $B(y_k,C\alpha_{n_k})$ and consists of $\approx \alpha_{n_k}^{s_j(1-1/\theta)}$ points that are $\gtrsim \alpha_{n_k}^{1/\theta}$-separated. Since these projected points lie on a line and are $\gtrsim \alpha_{n_k}^{1/\theta}$-separated, each interval of radius $(C\alpha_{n_k})^{1/\theta}$ contains at most a bounded number of them. Therefore
\begin{equation*}
  N_{(C\alpha_{n_k})^{1/\theta}}(P_V X \cap B(y_k,C\alpha_{n_k})) \ \gtrsim \ \alpha_{n_k}^{s_j(1-1/\theta)}
\end{equation*}
for infinitely many $k$, and the condition \eqref{alpha0} yields
\begin{equation*}
  \as P_V X \ \ge \ s_j
\end{equation*}
for all $V\in I_j$.

To prove the upper bound, fix again $V \in I_j$, set $r = R^{1/\theta}$, and consider $P_V X \cap B(x,R)$ for $x \in P_V X$ and $R \in (0,1)$. Because every direction in $\bigcup_{l=1}^\infty I_l$ lies in the fixed ball of radius $1/100$ about the first coordinate axis, the preliminary estimate on the separation of projected bricks from Section~\ref{sec:projection-two-values-upper} applies uniformly. Hence there exists $n_1 \in \N$ such that whenever $m > n \geq n_1$ are used indices, the projected bricks with indices $n$ and $m$ are separated by distance more than $2$: the proof in Section~\ref{sec:projection-two-values-upper} gives this for consecutive indices, and summing the positive projected increments only increases the gap. Therefore every ball $B(x,R)$ with $R \in (0,1)$ meets the projection of at most one such brick. The finitely many smaller-index bricks will be absorbed into the multiplicative constant at the end.

If the intersecting brick is $Q_{j,k,\ell}$ with index $n = n_{j,k,\ell}$, then its short axis is contained in $I_j$, so the same-family estimate from the proof of the even-brick case in Section~\ref{sec:projection-two-values-upper} gives
\begin{equation*}
  N_r(P_V X \cap B(x,R)) \ \lesssim \ \biggl(\frac{R}{r}\biggr)^{s_j+\delta_n}.
\end{equation*}
If the intersecting brick is $Q_{l,k',\ell'}$ with $l \neq j$ and index $n = n_{l,k',\ell'}$, then its short axis lies in $I_l$ and therefore makes an angle at least $d_j$ with $V$. In brick coordinates the projection is of the form $(u,v)\mapsto \eta u + \lambda v$ with $|\eta| \lesssim 1$ and $c_j \leq |\lambda| \lesssim 1$, where the constants depend only on $d_j$. Because \eqref{longgap} is uniform in the exponent $u \in (0,1]$, the projected contribution of the short side is eventually dominated by the minimal projected gap coming from the long coordinates. The odd-brick argument from Section~\ref{sec:projection-two-values-upper} therefore shows that
\begin{equation*}
  N_r(P_V X \cap B(x,R)) \ \lesssim_j \ \biggl(\frac{R}{r}\biggr)^{\delta_n}.
\end{equation*}
Fix $\varepsilon > 0$ and choose $n_0(j)$ so large that $\delta_n < \varepsilon$ for all $n \geq n_0(j)$. Then every brick with index at least $n_0(j)$ contributes at most
\begin{equation*}
  \lesssim_j \ \biggl(\frac{R}{r}\biggr)^{s_j+\varepsilon}
\end{equation*}
covering balls, while the finitely many smaller-index bricks can be absorbed into the multiplicative constant. Hence
\begin{equation*}
  \as P_V X \ \le \ s_j
\end{equation*}
finishing the proof of Theorem \ref{main1+}.

\section{Projections realising distinct full spectrum:~proof of Theorem \ref{maintheta}} \label{sec:projection-simultaneous-theta}

To prove Theorem \ref{maintheta}, we interlace, for a countable dense set of parameters $\theta$, the bricks used in the proof of Theorem \ref{main1}. Since the short-side exponents are constant, the covering estimates for the short and long coordinate sets from Section~\ref{sec:construction-of-the-set}, together with a case split according to whether the global parameter $\theta$ is above or below the parameter $\theta_i$ attached to a given brick, control the upper bound, while the dense set of parameters provides the lower bound at each $\theta$.

\subsection{Construction of the set}

Recall that the global sequences $(\alpha_n)_{n \in \N}$ and $(\beta_n)_{n \in \N}$ satisfy \eqref{alpha0} and \eqref{betasmall} for all $\theta \in (0,1)$. We keep the same disjoint open sets $I_s$ and $I_t$ from Section~\ref{sec:construction-of-the-set}, namely the balls of radius $1/100$ about the first and second coordinate axes. Let $\{\theta_i\}_{i = 1}^\infty$ be a countable dense subset of $(0,1)$.

For $u \in \{s,t\}$, $i \in \N$, and $n \in \N$, let $N_i(n,u)$ be chosen so that
\begin{equation*}
  2^{N_i(n,u)-1} \ \leq \ \alpha_n^{u(1-1/\theta_i)} \ \leq \ 2^{N_i(n,u)}.
\end{equation*}
Let $Q_i(n,u)$ denote the brick obtained by repeating the construction from Section~\ref{sec:construction-of-the-set} with $N = N_i(n,u)$. Thus the short coordinates of $Q_i(n,u)$ form a scaled copy of $\{\sum_{k = 1}^{N_i(n,u)} \pm 2^{-k/u}\}$ of diameter $\alpha_n$, the long coordinates form a scaled copy of $\{\sum_{k = 1}^{N_i(n,u)} \pm 2^{-k/\delta_n}\}$ of diameter $\beta_n$, and the set lies on the graph of a strictly increasing function. Since $\theta_i$ is fixed, the argument from Section~\ref{sec:construction-of-the-set} gives $N_i(n,u) = O_i(\log \log n)$ as $n \to \infty$, uniformly in $u \in \{s,t\}$. Hence $N_i(n,u)/\delta_n = o_i(\log n)$ and $2^{-N_i(n,u)/\delta_n} = n^{-o_i(1)}$, uniformly in $u \in \{s,t\}$. In particular,
\begin{equation*}
  \beta_n 2^{-N_i(n,u)/\delta_n} \ \to \ \infty
\end{equation*}
uniformly in $u \in \{s,t\}$. The same computation leading to \eqref{longgap} also yields
\begin{equation*}
  \frac{\alpha_n}{\beta_n 2^{-N_i(n,u)/\delta_n}} \ \to \ 0
\end{equation*}
as $n \to \infty$, uniformly in $u \in \{s,t\}$. Therefore, for each $i \in \N$, we may choose $k_0(i) \in \N$ so large that whenever $k \geq k_0(i)$, $n \in \{2^k,\ldots,2^{k+1}-1\}$, and $u \in \{s,t\}$, the separation argument from Section~\ref{sec:construction-of-the-set} implies that $Q_i(n,u)$ is $1$-separated and
\begin{equation*}
  \beta_n 2^{-N_i(n,u)/\delta_n} \ \geq \ 100\alpha_n.
\end{equation*}
The numerical constant is immaterial here; any sufficiently large absolute constant would suffice.
Choose pairwise disjoint infinite subsets $K_i^s$ and $K_i^t$ of $\{k_0(i),k_0(i)+1,\ldots\}$.

For each $i \in \N$ and each $k \in K_i^s$, let $\{V_{i,k,\ell}^s\}_{\ell = 1}^{M_{i,k}^s} \subset I_s$ be a maximal $2^{-k}$-separated subset of $I_s$. Since $\diam I_s = 1/50$, we have $M_{i,k}^s \leq \lfloor 2^k/50 \rfloor + 1 \leq 2^k$. Choose distinct indices
\begin{equation*}
  n_{i,k,\ell}^s \ \in \ \{2^k,\ldots,2^{k+1}-1\}
\end{equation*}
for $\ell \in \{1,\ldots,M_{i,k}^s\}$, choose $\phi_{i,k,\ell}^s \in [0,\pi)$ such that
\begin{equation*}
  V_{i,k,\ell}^s \ = \ \textup{span}(e^{i\phi_{i,k,\ell}^s}),
\end{equation*}
and define
\begin{equation*}
  Q_{i,k,\ell}^s \ = \ Q_i(n_{i,k,\ell}^s,s)e^{i\phi_{i,k,\ell}^s} + x_{n_{i,k,\ell}^s}.
\end{equation*}
Similarly, for each $i \in \N$ and each $k \in K_i^t$, let $\{V_{i,k,\ell}^t\}_{\ell = 1}^{M_{i,k}^t} \subset I_t$ be a maximal $2^{-k}$-separated subset of $I_t$. Since $\diam I_t = 1/50$, we have $M_{i,k}^t \leq \lfloor 2^k/50 \rfloor + 1 \leq 2^k$. Choose distinct indices
\begin{equation*}
  n_{i,k,\ell}^t \ \in \ \{2^k,\ldots,2^{k+1}-1\},
\end{equation*}
choose $\phi_{i,k,\ell}^t \in [0,\pi)$ such that
\begin{equation*}
  V_{i,k,\ell}^t \ = \ \textup{span}(e^{i\phi_{i,k,\ell}^t}),
\end{equation*}
and define
\begin{equation*}
  Q_{i,k,\ell}^t \ = \ Q_i(n_{i,k,\ell}^t,t)e^{i\phi_{i,k,\ell}^t} + x_{n_{i,k,\ell}^t}.
\end{equation*}
Finally, set
\begin{equation*}
  X \ = \ \bigcup_{i = 1}^\infty \biggl( \bigcup_{k \in K_i^s} \ \bigcup_{\ell = 1}^{M_{i,k}^s} Q_{i,k,\ell}^s \ \cup \ \bigcup_{k \in K_i^t} \ \bigcup_{\ell = 1}^{M_{i,k}^t} Q_{i,k,\ell}^t \biggr).
\end{equation*}
Because the sets $K_i^s$ and $K_i^t$ are pairwise disjoint and the dyadic blocks $\{2^k,\ldots,2^{k+1}-1\}$ are disjoint, each index $n$ is used for at most one brick. By the choice of $k_0(i)$, every used brick is $1$-separated. The separation argument from Section~\ref{sec:construction-of-the-set} also shows that distinct bricks are separated by distance at least $1$, so $X$ is uniformly discrete. Moreover, the maximal $2^{-k}$-separated subsets of $I_s$ and $I_t$ are non-empty, so whenever $k \in K_i^s$ or $k \in K_i^t$ there is at least one brick with index in $\{2^k,\ldots,2^{k+1}-1\}$. Since every family $K_i^s$ and $K_i^t$ is infinite and $|x_n| \to \infty$ with $n$, it follows that $X$ is unbounded. In particular, $X$ is locally finite and hence closed.

\subsection{Bounding the Assouad spectrum of projections}

Fix $\theta \in (0,1)$ and $V \in I_s$. Choose a sequence $(\theta_{i_m})_{m = 1}^\infty$ with $\theta_{i_m} \to \theta$ and $\theta_{i_m} \geq \theta$ for all $m$. For each $m$, choose $k_m \in K_{i_m}^s$ such that $k_m \to \infty$ as $m \to \infty$. Since $\{V_{i_m,k_m,\ell}^s\}_{\ell = 1}^{M_{i_m,k_m}^s}$ is a maximal $2^{-k_m}$-separated subset of $I_s$, there exists $\ell_m$ such that
\begin{equation*}
  \|V - V_{i_m,k_m,\ell_m}^s\| \ \lesssim \ 2^{-k_m} \ \lesssim \ \frac{1}{n_m},
\end{equation*}
where $n_m = n_{i_m,k_m,\ell_m}^s \in \{2^{k_m},\ldots,2^{k_m+1}-1\}$. Since $\theta_{i_m} \geq \theta$ and $0 < \alpha_{n_m} < 1$, we have $\alpha_{n_m}^{1/\theta_{i_m}} \geq \alpha_{n_m}^{1/\theta}$, and therefore \eqref{betasmall} gives $\beta_{n_m}/(n_m\alpha_{n_m}^{1/\theta_{i_m}}) \to 0$. Applying the lower-bound argument from Section~\ref{sec:projection-two-values-lower} to the selected brick $Q_{i_m,k_m,\ell_m}^s$, we obtain the same intra-brick estimate as before: there exist points $y_m \in P_V(Q_{i_m,k_m,\ell_m}^s) \subseteq P_V X$ and an absolute constant $C \geq 1$, independent of $m$, such that
\begin{equation*}
  P_V(Q_{i_m,k_m,\ell_m}^s) \ \subseteq \ B(y_m,C\alpha_{n_m})
\end{equation*}
for all sufficiently large $m$, and the set on the left consists of $\approx \alpha_{n_m}^{s(1-1/\theta_{i_m})}$ points which are $\gtrsim \alpha_{n_m}^{1/\theta_{i_m}}$-separated. Since $\theta_{i_m} \geq \theta$, the points are also $\gtrsim \alpha_{n_m}^{1/\theta}$-separated, and therefore, by monotonicity of the covering number in the set being covered,
\begin{equation*}
  N_{(C\alpha_{n_m})^{1/\theta}}(P_V X \cap B(y_m,C\alpha_{n_m})) \ \gtrsim \ \alpha_{n_m}^{s(1-1/\theta_{i_m})}.
\end{equation*}
Fix $\varepsilon > 0$. Since $\theta_{i_m} \to \theta$, we have $s(1-1/\theta_{i_m}) \to s(1-1/\theta)$, and hence
\begin{equation*}
  s(1-1/\theta_{i_m}) \ \leq \ (s-\varepsilon/2)(1-1/\theta)
\end{equation*}
for all sufficiently large $m$. Since $0 < \alpha_{n_m} < 1$, the preceding lower bound and the exponent inequality give
\begin{equation*}
  N_{(C\alpha_{n_m})^{1/\theta}}(P_V X \cap B(y_m,C\alpha_{n_m})) \ \gtrsim \ \alpha_{n_m}^{(s-\varepsilon/2)(1-1/\theta)} \ \approx \ \biggl(\frac{C\alpha_{n_m}}{(C\alpha_{n_m})^{1/\theta}}\biggr)^{s-\varepsilon/2}.
\end{equation*}
Because $C\alpha_{n_m} \to 0$, we also have
\begin{equation*}
  \biggl(\frac{C\alpha_{n_m}}{(C\alpha_{n_m})^{1/\theta}}\biggr)^{\varepsilon/2} \ \to \ \infty.
\end{equation*}
Therefore, given any $A > 0$, we have
\begin{equation*}
  N_{(C\alpha_{n_m})^{1/\theta}}(P_V X \cap B(y_m,C\alpha_{n_m})) \ \geq \ A \biggl(\frac{C\alpha_{n_m}}{(C\alpha_{n_m})^{1/\theta}}\biggr)^{s-\varepsilon}
\end{equation*}
for all sufficiently large $m$. Hence no constant works in the definition of $\as P_V X$ at the exponent $s-\varepsilon$, so $\as P_V X \geq s-\varepsilon$. Since $\varepsilon > 0$ was arbitrary, we conclude that $\as P_V X \geq s$ for all $\theta \in (0,1)$ and all $V \in I_s$. The same argument, with the roles of $s$ and $t$ interchanged, shows that $\as P_V X \geq t$ for all $\theta \in (0,1)$ and all $V \in I_t$.

To prove the upper bound, fix again $\theta \in (0,1)$ and $V \in I_s$, set $r = R^{1/\theta}$, and consider $P_V X \cap B(x,R)$ for $x \in P_V X$ and $R \in (0,1)$. Since $\beta_{n+1} \leq 2\beta_n$ for all $n$ and $\alpha_n/\beta_n \to 0$, the preliminary estimate on the separation of projected bricks from Section~\ref{sec:projection-two-values-upper} applies exactly as before. Hence there exists $n_0 \in \N$ such that any two projected bricks with used indices at least $n_0$ are separated by distance more than $2$, and therefore every ball $B(x,R)$ with $R \in (0,1)$ meets the projection of at most one brick with index at least $n_0$. The finitely many smaller-index bricks can be absorbed into the multiplicative constant.

Suppose first that the intersecting brick is $Q_{i,k,\ell}^s$ with index $n = n_{i,k,\ell}^s$, and identify $Q_i(n,s)$ with a subset of the $uv$-plane where the $u$-axis is the short side and the $v$-axis is the long side before rotation. Since $V$ and $V_{i,k,\ell}^s$ both lie in $I_s$, there exist absolute constants $0 < c_2 \leq C_2 < \infty$ and numbers $a_n,b_n,\lambda_n \in \R$ with $c_2 \leq |a_n| \leq C_2$ and $|\lambda_n| \leq C_2$ such that
\begin{equation*}
  P_V((u,v)e^{i\phi_{i,k,\ell}^s} + x_n) \ = \ P_V x_n + b_n + a_n(u + \lambda_n v)
\end{equation*}
for all $(u,v) \in Q_i(n,s)$. Replacing both $R$ and $r$ by comparable multiples changes covering numbers by at most an absolute constant, so it is enough to estimate the model map $(u,v) \mapsto u + \lambda_n v$. If $\alpha_n < R$, then the even-brick argument from Section~\ref{sec:projection-two-values-upper} applies with $Q_i(n,s)$ in place of $Q(2n,s)$, since that part of the proof uses only the affine model above and the covering estimates for the short and long coordinate sets, which are uniform in the truncation depth $N_i(n,s)$, and therefore gives
\begin{equation*}
  N_r(P_V X \cap B(x,R)) \ \lesssim \ \biggl(\frac{R}{r}\biggr)^{s+\delta_n}.
\end{equation*}
Assume from now on that $\alpha_n \geq R$. If $\theta_i > \theta$, then
\begin{equation*}
  r \ = \ R^{1/\theta} \ \leq \ \alpha_n^{1/\theta} \ < \ \alpha_n^{1/\theta_i},
\end{equation*}
so the relevant points can be covered individually. Since $Q_i(n,s)$ has $2^{N_i(n,s)}$ points, we obtain
\begin{equation*}
  N_r(P_V X \cap B(x,R)) \ \leq \ 2^{N_i(n,s)} \ \approx \ \alpha_n^{s(1-1/\theta_i)} \ \leq \ R^{s(1-1/\theta_i)}.
\end{equation*}
Writing
\begin{equation*}
  R^{s(1-1/\theta_i)} \ = \ \biggl(\frac{R}{r}\biggr)^{s\theta(1-\theta_i)/(\theta_i(1-\theta))}
\end{equation*}
and using $\theta < \theta_i$, we conclude that
\begin{equation*}
  N_r(P_V X \cap B(x,R)) \ \lesssim \ \biggl(\frac{R}{r}\biggr)^s.
\end{equation*}
Suppose instead that $\theta_i \leq \theta$. If $\lambda_n = 0$, then the relevant first coordinates lie in an interval of length $\lesssim R$, and the covering estimate from Section~\ref{sec:construction-of-the-set} gives
\begin{equation*}
  N_r(P_V X \cap B(x,R)) \ \lesssim \ \biggl(\frac{R}{r}\biggr)^s.
\end{equation*}
Assume now that $\lambda_n \neq 0$. The relevant points of $Q_i(n,s)$ satisfy
\begin{equation*}
  |u + \lambda_n v - c| \ \leq \ CR
\end{equation*}
for some $c \in \R$. Since $0 \leq u \leq \alpha_n$ and $\alpha_n \geq R$, it follows that $\lambda_n v$ lies in an interval of length $\lesssim \alpha_n$. Hence the corresponding long coordinates lie in an interval $J$ of length $\lesssim \alpha_n/|\lambda_n|$. We cover the intersection of the long coordinates with $J$ by intervals of length $\min\{r/|\lambda_n|,\beta_n\}$. Since the long coordinates of $Q_i(n,s)$ form a scaled copy of a set of the form $\{\sum_{k = 1}^{N_i(n,s)} \pm 2^{-k/\delta_n}\}$ of diameter $\beta_n$, the covering estimate from Section~\ref{sec:construction-of-the-set} gives
\begin{equation*}
  N_v \ \lesssim \ \biggl(\frac{\min\{\alpha_n/|\lambda_n|,\beta_n\}}{\min\{r/|\lambda_n|,\beta_n\}}\biggr)^{\delta_n}.
\end{equation*}
If $|\lambda_n|\beta_n \leq r$, then both minima equal $\beta_n$ and $N_v \lesssim 1$. If $r < |\lambda_n|\beta_n \leq \alpha_n$, then the numerator equals $\beta_n$ and the denominator equals $r/|\lambda_n|$, so
\begin{equation*}
  N_v \ \lesssim \ \biggl(\frac{|\lambda_n|\beta_n}{r}\biggr)^{\delta_n} \ \leq \ \biggl(\frac{\alpha_n}{r}\biggr)^{\delta_n}.
\end{equation*}
If $|\lambda_n|\beta_n > \alpha_n$, then the numerator equals $\alpha_n/|\lambda_n|$ and the denominator equals $r/|\lambda_n|$, so
\begin{equation*}
  N_v \ \lesssim \ \biggl(\frac{\alpha_n}{r}\biggr)^{\delta_n}.
\end{equation*}
Thus in every case
\begin{equation*}
  N_v \ \lesssim \ \biggl(\frac{\alpha_n}{r}\biggr)^{\delta_n}.
\end{equation*}
For each interval $W$ in this cover, the variation of $\lambda_n v$ over $W$ is at most $r \leq R$. Therefore the corresponding first coordinates lie in an interval of length $\lesssim R$, and the covering estimate for the short coordinate set gives at most
\begin{equation*}
  \lesssim \ \biggl(\frac{R}{r}\biggr)^s
\end{equation*}
balls. Consequently,
\begin{equation*}
  N_r(P_V X \cap B(x,R)) \ \lesssim \ \biggl(\frac{\alpha_n}{r}\biggr)^{\delta_n}\biggl(\frac{R}{r}\biggr)^s.
\end{equation*}
If the intersecting brick is $Q_{i,k,\ell}^t$ with index $n = n_{i,k,\ell}^t$, identify $Q_i(n,t)$ with a subset of the $uv$-plane where the $u$-axis is the short side and the $v$-axis is the long side before rotation. Since $V \in I_s$ while $V_{i,k,\ell}^t \in I_t$, there exist absolute constants $C_3 < \infty$ and $c_3 > 0$ such that the projection onto $V$ is of the form
\begin{equation*}
  (u,v) \ \mapsto \ \eta_n u + \lambda_n v
\end{equation*}
with $|\eta_n| \leq C_3$ and $c_3 \leq |\lambda_n| \leq C_3$ for all $n$. Since $k \in K_i^t \subseteq \{k_0(i),k_0(i)+1,\ldots\}$, the choice of $k_0(i)$ and the estimate $\beta_n 2^{-N_i(n,t)/\delta_n} \geq 100\alpha_n$ show that the projected contribution of the short coordinate is dominated by the projected long-coordinate gap. The numerical constant $100$ is chosen large enough, relative to the uniform bounds on $\eta_n$ and $\lambda_n$, that for all distinct pairs of points $(u_i,v_i),(u_j,v_j) \in Q_i(n,t)$ we have
\begin{equation*}
  |\eta_n||u_j-u_i| \ \leq \ \frac12 |\lambda_n||v_j-v_i|.
\end{equation*}
Hence the projected brick is bi-Lipschitz to a scaled copy of the long coordinate set. Therefore every interval of length $R$ in the projection corresponds to an interval of long coordinates of length $\lesssim R$, and the odd-brick argument from Section~\ref{sec:projection-two-values-upper} yields
\begin{equation*}
  N_r(P_V X \cap B(x,R)) \ \lesssim \ \biggl(\frac{R}{r}\biggr)^{\delta_n}.
\end{equation*}
Fix $\varepsilon > 0$ and choose $n_1 \geq n_0$ so large that $\delta_n < \varepsilon(1-\theta)$ for all $n \geq n_1$; recall that $R/r = R^{1-1/\theta} \geq 1$. Suppose first that the intersecting brick is a same-family $s$-brick with index at least $n_1$. In the case $\alpha_n < R$, the estimate above already gives
\begin{equation*}
  N_r(P_V X \cap B(x,R)) \ \lesssim \ \biggl(\frac{R}{r}\biggr)^{s+\delta_n} \ \leq \ \biggl(\frac{R}{r}\biggr)^{s+\varepsilon}.
\end{equation*}
If $\alpha_n \geq R$ and $\theta_i > \theta$, then the point-count estimate above gives
\begin{equation*}
  N_r(P_V X \cap B(x,R)) \ \lesssim \ \biggl(\frac{R}{r}\biggr)^s \ \leq \ \biggl(\frac{R}{r}\biggr)^{s+\varepsilon}.
\end{equation*}
If $\alpha_n \geq R$ and $\theta_i \leq \theta$, then $\alpha_n < 1$, so the extra long-coordinate factor in the bound above satisfies
\begin{equation*}
  \biggl(\frac{\alpha_n}{r}\biggr)^{\delta_n} \ \leq \ r^{-\delta_n} \ = \ R^{-\delta_n/\theta} \ \leq \ \biggl(\frac{R}{r}\biggr)^{\varepsilon},
\end{equation*}
and therefore $N_r(P_V X \cap B(x,R)) \lesssim (R/r)^{s+\varepsilon}$ in this case as well. Thus, whenever the intersecting brick is a same-family $s$-brick with index at least $n_1$,
\begin{equation*}
  N_r(P_V X \cap B(x,R)) \ \lesssim \ \biggl(\frac{R}{r}\biggr)^{s+\varepsilon}.
\end{equation*}
Suppose instead that the intersecting brick is an opposite-family $t$-brick with index at least $n_1$. Then
\begin{equation*}
  N_r(P_V X \cap B(x,R)) \ \lesssim \ \biggl(\frac{R}{r}\biggr)^{\delta_n} \ \leq \ \biggl(\frac{R}{r}\biggr)^{\varepsilon}.
\end{equation*}
The finitely many smaller-index bricks are absorbed into the constant. Hence $\as P_V X \leq s$ for all $\theta \in (0,1)$ and all $V \in I_s$. By the same argument, with the roles of $s$ and $t$ interchanged, we obtain $\as P_V X \leq t$ for all $\theta \in (0,1)$ and all $V \in I_t$. This completes the proof of Theorem \ref{maintheta}.

\section*{Acknowledgements}

We thank Stuart Burrell for helpful discussions about the topic and for giving us access to some typed notes from 2021 concerning Theorem~\ref{main2}.


\begin{thebibliography}{Cana}

  \bibitem[AHRS21]{anderson} T. Anderson, K. Hughes, J. Roos and A. Seeger. 	$L^p\to L^q$ bounds for spherical maximal operators,
  \emph{Math. Z.}, 297, (2021), 1057--1074.



  \bibitem[BCP24+]{holder}
  A. Baraviera, M. Carvalho, and G. Pessil.
  Metric mean dimension, Hölder regularity and Assouad spectrum, preprint: arXiv:2407.15774

  \bibitem[BRS24+]{beltran}
  D.  Beltran, J.  Roos, and A. Seeger.
  Spherical maximal operators with fractal sets of dilations on radial functions, preprint: arXiv:2412.09390

  \bibitem[CG24]{stathisquasi}
  E. K. Chrontsios Garitsis.
  Quasiregular distortion of dimensions,
  \emph{Conform. Geom. Dyn.}, {\bf 28}, (2024), 165--175.

  \bibitem[CGT24]{stathis} E. K. Chrontsios Garitsis  and  S. Troscheit. Minkowski weak embedding theorem, \emph{Houston J. Math.}, {\bf 50}, (2024), 259--273.

  \bibitem[CGT23]{stathistyson}
  E. K. Chrontsios  Garitsis, and J. Tyson. Quasiconformal distortion of the Assouad spectrum and classification of polynomial spirals,
  \emph{Bull. Lond. Math. Soc.}, {\bf 55}, (2023), 282--307.

  \bibitem[F89]{Falconer1989}
  K. J. Falconer. Dimensions and measures of quasi self-similar sets,
  \emph{Proc. Amer. Math. Soc.}, {\bf 106}, (1989), 543--554.

  \bibitem[F14]{Fal03}K. J. Falconer. \emph{Fractal Geometry: Mathematical Foundations and Applications}. John Wiley and Sons, Hoboken, NJ, 3rd. ed., (2014).

  \bibitem[F20]{falconerprofile2}
  {K. J. Falconer}, A capacity approach to box and packing dimensions of projections and other images,  \emph{Analysis, Probability and Mathematical Physics on Fractals}, World Scientific, (2020).

  \bibitem[F21]{falconerprofile}
  {K. J. Falconer}, A capacity approach to box and packing dimensions of projections of sets and exceptional directions, {\em J. Fractal Geom.},  {\bf 8}, (2021), 1--26.

  \bibitem[F25+]{falconer25}K. J. Falconer. Seventy Years of Fractal Projections, \emph{Fractal Geometry and Stochastics VII}, Springer, (2026), available at arXiv:2602.22002.

  \bibitem[FFJ15]{FFJ15}K. J. Falconer, J. M. Fraser and X. Jin. Sixty Years of Fractal Projections. \emph{Fractal Geometry and Stochastics V}, Springer International Publishing, (2015) 3--25.

  \bibitem[FH96]{falconerhowroyd}
  K. J. Falconer and J. D. Howroyd.
  Projection theorems for box and packing dimensions,
  \emph{Math. Proc. Cambridge Philos. Soc.}, {\bf119}, (1996), 287--295.

  \bibitem[FH97]{falconerhowroyd2}
  {K. J. Falconer \and J. D. Howroyd}, Packing dimensions of projections and dimension profiles, {\em Math. Proc. Cambridge Philos. Soc.}, 121:269-286 (1997).

  \bibitem[F18]{fraserisrael}J. M. Fraser.
  Distance sets, orthogonal projections, and passing to weak tangents, {\it Israel J. Math.}, {\bf 226}, (2018), 851--875.

  \bibitem[Fr21a]{jon:book}
  J. M. Fraser. {\em Assouad Dimension and Fractal Geometry},
  Cambridge University Press,
  Tracts in Mathematics Series, 2021.

  \bibitem[FHHTY19]{canadian}
  J.~M. Fraser, K. E. Hare, K. G. Hare, S. Troscheit and H. Yu.
  The Assouad spectrum and the quasi-Assouad dimension: a tale of two spectra,
  \emph{Ann. Acad. Sci. Fenn. Math.}, {\bf 44}, (2019), 379--387.
  
    \bibitem[FHOR15]{fhor}
  J.~M. Fraser, A. M. Henderson, E. J. Olson and J. C. Robinson.  
  On the Assouad dimension of self-similar sets with overlaps,
        \emph{Adv. Math.}, {\bf 273}, (2015), 188--214.

  \bibitem[FK20]{antti}
  J.~M. Fraser and A. K\"aenm\"aki.
  Attainable values for the Assouad dimension of projections, \emph{Proc. Amer. Math. Soc.}, {\bf 148}, (2020), 3393--3405.

  \bibitem[FO17]{FO17}J. M. Fraser and T. Orponen. The Assouad dimensions of projections of planar sets. \emph{Proc. Lond. Math. Soc.}, \textbf{114}, (2017), 374--398.

  \bibitem[FS24]{bullams}
  J. M. Fraser and L. Stuart. A new perspective on the Sullivan dictionary via Assouad type dimensions and spectra. {\it
  Bull. Amer. Math. Soc.}, {\bf 61}, (2024), 103--118.

  \bibitem[FY18]{assouadspectrum}J. M. Fraser and H. Yu. New dimension spectra: Finer information on scaling and homogeneity. \emph{Adv. Math.}, \textbf{329}, (2018), 273--328.

  \bibitem[H01]{howroyd}
  J. D. Howroyd. Box and packing dimensions of projections and dimension profiles,
  \emph{Math. Proc. Cambridge Philos. Soc.}, {\bf 130}, (2001), 135--160.

  \bibitem[K68]{Kau68}R. Kaufman. On Hausdorff dimension of projections. \emph{Mathematika}, \textbf{15}, (1968), 153--155.

  \bibitem[KR23+]{KR23}
  A. K\"aenm\"aki and A. Rutar.
  Tangents of invariant sets, preprint: arXiv:2309.11971

  \bibitem[LX16]{quasiassouad}
  F. L\"u and L.-F. Xi.
  Quasi-Assouad dimension of fractals,
  {\em J. Fractal Geom.}, {\bf 3}, (2016), 187--215.

  \bibitem[MT10]{mackaytyson}
  J. M. Mackay and J. T. Tyson.
  \emph{Conformal dimension. Theory and application},
  University Lecture Series, 54. American Mathematical Society, Providence, RI, 2010.

  \bibitem[M54]{Mar54}J. Marstrand. Some Fundamental Geometrical Properties of Plane Sets of Fractional Dimensions. \emph{Proc. Lond. Math. Soc. (3)}, {\bf 4}, (1954), 257--302.

  \bibitem[M14]{Mat14}P. Mattila. Recent Progress on Dimensions of Projections. \emph{Geometry and Analysis of Fractals}, (2014), 283--301.

  \bibitem[M75]{Mat75}P. Mattila. Hausdorff dimension, orthogonal projections and intersections with planes. \emph{Ann. Acad. Sci. Fenn. Ser. A I Math.}, \textbf{1}(2), (1975), 227--244.

  \bibitem[O21]{orponenassouad} T. Orponen. On the Assouad dimension of projections. \emph{Proc. Lond. Math. Soc. (3)}, {\bf 122}, (2021), 317--351.

  \bibitem[R11]{robinson}
  J. C. Robinson.
  {\em Dimensions, Embeddings, and Attractors},
  Cambridge University Press, 2011.

  \bibitem[RS23]{roos} J. Roos and A. Seeger. Spherical maximal functions and fractal dimensions of dilation sets, \emph{Amer. J. Math.}, {\bf 145}, (2023), 1077--1110.

  \bibitem[R24]{specclass}
  A. Rutar.
  Attainable forms of Assouad spectra, \emph{Indiana Univ. Math. J.}, {\bf 73}, (2024), 1331--1356

  \bibitem[S15]{Shm15}P. Shmerkin. Projections of Self-Similar and Related Fractals: A Survey of Recent Developments. In: \emph{Fractal Geometry and Stochastics V}, Springer International Publishing, 53--74 (2015).

  \bibitem[WL25+]{wangli}
  L. Wang and W. Li.
  The generalized upper box dimension, preprint: arXiv:2510.00521

  \bibitem[W25+]{wu}
  M. Wu.
  Projection theorems with countably many exceptions and applications to the exact overlaps conjecture, preprint: arXiv:2503.21923

\end{thebibliography}
\end{document}